\documentclass{ws-m3as}
\usepackage{algpseudocode}
\usepackage{algorithm}
\usepackage{booktabs}

\begin{document}

\markboth{S. Wang, C. Liu, and W. Hao}{Learning Dynamical Systems via Energetic Variational Optimal Transport}

%
\catchline{}{}{}{}{}
%
\title{Learning the Energy Landscapes of Dynamical Systems via Energetic Variational Optimal Transport under Data Quantity--Quality Trade-offs}

\author{Shun Wang}

\address{Department of Mathematics, Penn State University, University Park, Pennsylvania, United States of America}

\author{Chun Liu\footnote{Corresponding author}}

\address{Department of Applied Mathematics, Illinois Institute of Technology, Chicago, United States of America}

\author{Wenrui Hao\footnote{Corresponding author}}

\address{Department of Mathematics, Penn State University, University Park, Pennsylvania, United States of America}

\maketitle


\begin{abstract}
Dynamic optimal transport unifies optimal transport, fluid mechanics, and gradient-flow theory within a continuous dynamical framework, offering a geometry-aware language that underpins a broad range of applications across physics, biology, and machine learning. However, conventional formulations cast it as a constrained optimization problem that must explicitly satisfy the continuity equation, which restricts modeling flexibility and hinders the reconstruction of the underlying dynamics directly from data. We propose the energetic variational method for dynamic optimal transport (EVMDOT), which reformulates dynamic optimal transport within an energetic variational framework by combining the flow map, the least action principle, and the maximum dissipation principle. The flow map recasts the constrained problem as an unconstrained one by automatically enforcing the continuity equation, while the balance between the conservative and dissipative forces, obtained by taking variations of the free energy and the dissipation, determines the velocity field. Applied to the Fokker--Planck equation, the EVMDOT reconstructs both the energy landscape and the Waddington landscape directly from time-series density data, with the potential parameterized by a multilayer neural network. Through numerical experiments, we reveal that the EVMDOT achieves an intrinsic balance between data quantity and data quality: a sufficient data quantity compensates for limited data quality, making the reconstruction robust to the choice of the observation window. We further apply the EVMDOT to the Alzheimer's Disease Neuroimaging Initiative (ADNI) dataset to infer the potential landscape of amyloid-$\beta$ and tau, revealing two wells corresponding to the cognitively normal and Alzheimer's disease stages and capturing the key transition pathway between them. The EVMDOT  establishes an extensible, data-driven framework whose energetic variational structure can be readily adapted to new systems, opening broad potential for quantifying disease progression and uncovering the dynamics that govern complex biological processes.
\end{abstract}


\section{Introduction}
Dynamic optimal transport (DOT) extends the classical optimal transport problem by modeling the evolution of probability distributions over time, providing a powerful framework for describing flows, resource allocation, and system dynamics across diverse fields, including physics, biology, control theory, and machine learning. The foundational Benamou--Brenier formulation \cite{benamou2000computational,benamou2002monge} recasts the static Wasserstein distance as a convex optimization problem constrained by the continuity equation, thereby enabling dynamic interpretations and efficient computational methods \cite{haasler2024scalable,van2023thermodynamic,gangbo2019unnormalized}. Recent advances have unified DOT with thermodynamic principles \cite{papadakis2014optimal} and extended it to networks \cite{tong2020trajectorynet}, discrete surfaces \cite{burger2023dynamic}, nonlinear control systems \cite{lavenant2018dynamical,cang2025synchronized}, and high-dimensional data via deep learning \cite{bunne2024optimal}. Its applications range from modeling cellular dynamics in single-cell omics \cite{schiebinger2019optimal,cang2023screening} to parameter identification in chaotic systems \cite{ghoussoub2018optimal}. Methodological innovations include scalable algorithms for high-dimensional problems \cite{bunne2024optimal,schiebinger2019optimal} and extensions to unbalanced or constrained transport settings \cite{chizat2018interpolating,chen2021optimal,sun2026variational,zhang2025learning}.

To make this concrete, we first recall the two equivalent formulations of dynamic optimal transport. In the static Monge formulation, the squared Wasserstein-2 distance between two probability density functions $p_1(X)$ and $p_2(x)$ is defined as
\begin{equation}
\label{eq:1}
W_2^2\big(p_1(X), p_2(x)\big) = \inf_{M} \int_{\Omega} \|X - M(X)\|^2 \, p_1(X)\, dX,
\end{equation}
where $M:X \rightarrow x$ is a smooth one-to-one map from $\Omega$ onto itself, satisfying the mass-preserving condition $\int_{x \in A} p_2(x)\, dx = \int_{M(X) \in A} p_1(X)\, dX$ for every bounded subset $A \subset \Omega$. While the Monge formulation seeks an optimal static map $M$, it is generally nonconvex and offers no direct description of how mass is transported over time. The Benamou--Brenier formulation \cite{benamou2000computational,benamou2002monge} overcomes these difficulties by recasting the static Monge problem as an equivalent dynamic optimization over a time-dependent density $\rho(x,t)$ and velocity field $u(x,t)$:
\begin{equation}
\label{eq:2}
W_2^2(\rho, u) = \inf_{\rho,\, u} \; T \int_0^T \int_{\Omega} \rho(x,t)\, \|u(x,t)\|^2 \, dx \, dt,
\end{equation}
\[
\text{s.t.} \quad
\begin{cases}
\partial_t \rho + \nabla \cdot (\rho u) = 0, \\
\rho(x,0) = p_1(X), \quad \rho(x,T) = p_2(x).
\end{cases}
\]
Here the continuity equation enforces mass conservation along the transport path, while the boundary conditions fix the initial and terminal densities. The optimal velocity field traces the Wasserstein geodesic connecting $p_1$ and $p_2$, transforming the static transport map into a continuous flow of probability mass and thereby endowing dynamic optimal transport with both a clear physical interpretation and a route to efficient computation.

Despite its elegant variational structure, the Benamou--Brenier formulation poses a constrained optimization problem in which the continuity equation $\partial_t \rho + \nabla \cdot (\rho u) = 0$ must be satisfied throughout the time interval. Enforcing this constraint explicitly is computationally demanding, as it couples the density $\rho(x,t)$ and the velocity field $\frac{dx}{dt}=u(x,t)$  at every point in space and time, and typically requires spatial discretization of the density that scales poorly with dimension. These difficulties limit the flexibility of the formulation and hinder its direct application to high-dimensional, data-driven settings. To alleviate these issues, continuous normalizing flows (CNFs) \cite{sha2024reconstructing,tong2020trajectorynet} parameterize the velocity field with a neural network and evolve the density as $\frac{d\ln\rho}{dt}=-\nabla\cdot u$ through a neural ordinary differential equation (neural ODE) \cite{chen2018neural}. By integrating the associated instantaneous change-of-variables formula along the flow, CNFs satisfy the continuity equation implicitly and thus avoid its explicit enforcement. However, evaluating the instantaneous change in density requires computing the Jacobian trace (i.e. the divergence) of the velocity field $u(x,t)$ at each integration step, which is computationally expensive in high dimensions. Moreover, parameterizing $u(x,t)$ with a generic neural network in place of an explicit energy gradient introduces substantial instability into the learned dynamics, and the theoretical guarantees established under idealized assumptions often fail to hold on real, noisy data. As a further consequence, the absence of an underlying energetic structure prevents the recovery of the potential landscape that governs the system.

By introducing a scalar Lagrange multiplier $\phi(x,t)$, which plays the role of a scalar potential function, the dual problem of the Benamou--Brenier formulation can be written as a Hamilton--Jacobi equation in $\phi$ \cite{su1999discontinuous,benamou2000computational,benamou2002monge} (see Appendix for details):
\begin{equation}
\label{eq:3}
\sup_{\phi} \; T \int_{\Omega} \big[\phi(x,T)\rho(x,T) - \phi(x,0)\rho(x,0)\big]\, dx,
\end{equation}
\[
\text{s.t.} \quad
\partial_t \phi + \frac{1}{2}|\nabla\phi|^2 = 0.
\]
The fact that the dual of dynamic optimal transport takes the form of a Hamilton--Jacobi equation reveals a deep connection between optimal transport and classical mechanics, and motivates us to reexamine dynamic optimal transport from the perspective of variational mechanics. 

Building on this insight, we propose the energetic variational method \cite{xu2014energetic,lu2024learning,giga2017variational} for dynamic optimal transport (EVMDOT), which reformulates dynamic optimal transport within an energetic variational framework by combining the flow map, the least action principle (LAP), and the maximum dissipation principle (MDP). In this formulation, the flow map recasts the constrained problem as an unconstrained one by automatically enforcing the continuity equation, while the force balance between the conservative and dissipative forces, derived respectively from the free energy and the dissipation, determines the velocity field and recovers the Hamilton--Jacobi structure from the perspective of variational mechanics. We demonstrate the interpretability and broad applicability of EVMDOT through a series of applications, ranging from reconstructing energy and Waddington landscapes  for the Fokker--Planck equation (See Algorithm \ref{alg:1}) to inferring the potential landscape of amyloid-$\beta$ and tau from ADNI datasets. These applications establish EVMDOT as an interpretable, data-driven framework.

\section{Methods}
\subsection{Flow Map}
For a given velocity field $u(x,t)$, the corresponding flow map $x(X,t)$ is defined by
\begin{equation}
x_t = u(x,t), x(X,0) = X.
\end{equation}
Here, $x(X,t)\in\Omega_x^t$ denotes the Eulerian coordinates, and $X\in\Omega_X^0$ denotes the Lagrangian coordinates, i.e. the initial configuration. To characterize the evolution of structures or patterns (the configuration), the deformation gradient \cite{giga2017variational} is defined as
\begin{equation}
F(X,t) = \frac{\partial x(X,t)}{\partial X}
\end{equation}
In Eulerian coordinates, for a fluid of density $\rho(x,t)$, the relation $\rho(x,t)=\rho(X,0)/\mathrm{det} F$ (equivalently, $\int_{\Omega^t_x} \rho(x,t)dx = \int_{\Omega^0_X} \rho(X,0)dX$) yields the local mass conservation law
\begin{equation}
\rho_t+\nabla_x\cdot(\rho u) = 0,
\end{equation}
whose derivation proceeds as follows:
\begin{align}
&\frac{d\rho}{dt} = \partial_t\rho+(u\cdot\nabla)\rho =-\rho(X,0)\frac{1}{(\mathrm{det}F)^2}\mathrm{det}F\;\mathrm{tr}(F^{-1}\partial_tF) = -\rho(x,t) \mathrm{tr}(F^{-1}\partial_tF)  \\
& \mathrm{tr}(F^{-1}\partial_tF)=\mathrm{tr}(\frac{\partial X}{\partial x}\frac{\partial u}{\partial X})=\mathrm{tr}(\frac{\partial u}{\partial x})=\nabla\cdot u \\
&\partial_t\rho+(u\cdot\nabla)\rho + \rho\nabla\cdot u = \rho_t+\nabla\cdot(\rho u) =  0
\end{align}

\subsection{Energetic Variational Method for Dynamic Optimal Transport (EVMDOT)}
By combining the least action principle (LAP), the maximum dissipation principle (MDP), and Newton's third law, the Hamilton--Jacobi formulation of dynamic optimal transport can be derived. First, the action functional $\mathcal{A}$ is defined as
\begin{equation}
\mathcal{A} (\phi,\rho)= \int_{0}^{T}\int_{\Omega} \phi(x,t)\rho(x,t)dxdt,
\end{equation}
and the dissipation functional $\mathcal{D}$ as
\begin{equation}
\mathcal{D} (u,\rho)= \frac{1}{2}\int_{0}^{T}\int_{\Omega} \|u(x,t)\|^2\rho(x,t)dxdt.
\end{equation}
The energy dissipation $\mathcal{D}(u,\rho)$ corresponds to Darcy's law in fluid dynamics (friction relative to the resting medium), so that dynamic optimal transport is equivalent to the MDP ($\inf_{u,\rho} \; \mathcal{D}(u,\rho)$). The equivalent formulation of the mass conservation law reads
\begin{equation}
\rho(x,t) = \rho(X,0)/\mathrm{det} F(X,t).
\end{equation}
Applying the LAP in the space $H=L^2(0,T;\tilde{H})$ equipped with an inner product $\langle \cdot\rangle$, the force $\mathrm{force}_{\mathcal{A}}$ is written as
\begin{equation}
\mathrm{force}_{\mathcal{A}} = H_{-}\frac{\delta \mathcal{A}}{\delta x},
\end{equation}
which gives
\begin{equation}
\mathrm{force}_{\mathcal{A}} = \rho(x,t)\nabla\phi(x,t).
\end{equation}
The derivation proceeds as follows:
\begin{align}
\delta \mathcal{A} &=\delta\int_{0}^{T}\int_{\Omega} \phi(x,t)\rho(x,t)dxdt \nonumber \\
&=\delta\int_{0}^{T}\int_{\Omega} \phi(x,t)(\rho(X,0)/\mathrm{det} F(X,t)) \mathrm{det} F(X,t)dXdt \nonumber \\
&= \int_{0}^{T}\int_{\Omega} (\nabla\phi(x,t)\cdot\delta x)\rho(X,0)dXdt =  \int_{0}^{T}\int_{\Omega} \rho(x,t)\nabla\phi(x,t)\cdot\delta xdxdt.
\end{align}
The MDP states that the dissipative force is obtained by minimizing the dissipation functional $\mathcal{D}$ with respect to the aforementioned ``rate''. Accordingly, the dissipative force is derived as
\begin{equation}
\mathrm{force}_{dissipative} = \tilde{H}_{-}\frac{\delta \mathcal{D}}{\delta x_t}.
\end{equation}
and is expressed as
\begin{equation}
\mathrm{force}_\mathrm{dissipative} = \rho(x,t)u(x,t).
\end{equation}
The derivation proceeds as follows:
\begin{align}
\delta \mathcal{D} &= \frac{1}{2}\delta \int_{0}^{T}\int_{\Omega} \|u(x,t)\|^2\rho(x,t)dxdt \nonumber \\
&= \frac{1}{2}\delta\int_{0}^{T}\int_{\Omega} \|u(x,t)\|^2(\rho(X,0)/\mathrm{det} F(X,t)) \mathrm{det} F(X,t)dXdt \nonumber \\
&= \int_{0}^{T}\int_{\Omega} \rho(X,0)u(x,t)\cdot\delta udXdt= \int_{0}^{T}\int_{\Omega} \rho(x,t)u(x,t)\cdot\delta udxdt;
\end{align}
By Newton's third law ($\mathrm{force}_{\mathcal{A}} =\mathrm{force}_\mathrm{dissipative}$), the velocity $u(x,t)$ equals the gradient of the action potential:
\begin{equation}
u(x,t) = \nabla\phi(x,t)
\end{equation}
According to the Hamilton--Jacobi equation, the action $\phi(x,t)$ satisfies
\begin{equation}
\partial_t\phi(x,t)+H(x,p,t) = 0
\end{equation}
where the Hamiltonian $H(x,p,t)$ obeys
\begin{equation}
x_t = \frac{\partial H}{\partial p},~ p = \frac{\partial \phi}{\partial x}.
\end{equation}
Substituting $u(x,t) = \nabla\phi(x,t)$ into the equations above, the Hamiltonian is obtained as
\begin{equation}
\frac{\partial H}{\partial p} = p, p = \nabla\phi(x,t),
H(x,p,t) = \frac{1}{2} |\nabla\phi(x,t)|^2.
\end{equation}
Hence, the Hamilton--Jacobi equation for dynamic optimal transport reads
\begin{equation}
\partial_t\phi(x,t)+\frac{1}{2} |\nabla\phi(x,t)|^2 = 0.
\end{equation}
We thereby recover the Hamilton--Jacobi equation of optimal transport from mechanics.

Integrating the flow map and the mechanical formulation, dynamic optimal transport can be recast into the following system:
\begin{equation}
\quad
\begin{cases}
\rho(x,t) = \rho(X,0)/\mathrm{det} F(X,t), \\
\delta_{x} \mathcal{A}=\delta_{x_t} \mathcal{D}.
\end{cases}
\end{equation}
Here, $x$ denotes the trajectory in Lagrangian coordinates. The action functional $\mathcal{A}$ is rewritten as
\[
\mathcal{A}(x) = \int_{0}^{T} (\mathcal{K}(x)-\mathcal{F}(x))dt 
\]
where $\mathcal{K}(x)$ is the kinetic energy and $\mathcal{F}(x)$ is the Helmholtz free energy. Applying the LAP in the space $H=L^2(0,T;\tilde{H})$ equipped with an inner product $\langle \cdot\rangle$, the variations $\delta_x \mathcal{K}$ and $\delta_x \mathcal{F}$ yield
 \begin{align}
  &\mathrm{force}_{intertial} = H_{-}\frac{\delta \int_{0}^{T} \mathcal{K}dt}{\delta x}; \nonumber \\
 &\mathrm{force}_{conservative} = H_{-}\frac{\delta \int_{0}^{T} \mathcal{F}dt}{\delta x}.
 \end{align}
Once all forces have been derived, the force balance follows from Newton's third law:
\begin{equation}
\mathrm{force}_{intertial} = \mathrm{force}_{conservative}+\mathrm{force}_{dissipative}
\end{equation}
Since the kinetic energy is typically assumed to vanish, $\mathcal{K}=0$, the force balance $(\delta_x \mathcal{A}=\delta_{x_t} \mathcal{D})$ reduces to
\begin{equation}
\frac{\delta \mathcal{F}}{\delta x}+\frac{\delta \mathcal{D}}{\delta x_t}=0
\end{equation}
When $\mathcal{F}[\rho]=\int_\Omega \omega(\rho)dx$ and $\mathcal{D} (u,\rho)= \frac{1}{2}\int_{0}^{T}\int_{\Omega} \|u(x,t)\|^2\rho(x,t)dxdt$, the velocity field $u(x,t)$ is derived as
\begin{align}
&L^2_{x,t_{-}} \frac{\delta \mathcal{D}}{\delta x_t} = \rho u,\quad
L^2_{x,t_{-}} \frac{\delta \mathcal{F}}{\delta x} = \nabla_x[\omega_{\rho}(\rho)\rho-\omega(\rho)] \nonumber \\
&u(x,t) = -\frac{1}{\rho} \nabla_x[\omega_{\rho}(\rho)\rho-\omega(\rho)]    
\end{align}
The derivation proceeds as follows:
\begin{align}
\label{eq:4}
\delta \int_0^T \mathcal{F} \, dt 
&= \delta \int_0^T \int_\Omega \omega\left( \frac{\rho_0}{\det F} \right) \det F \, dX dt \nonumber \\
&= \int_0^T \int_\Omega \left[ -\omega_\rho \left( \frac{\rho_0}{\det F} \right) \frac{\rho_0}{\det F} + \omega\left( \frac{\rho_0}{\det F} \right) \right] \det F \, \text{tr} \left( F^{-1} \frac{\partial \delta x}{\partial X} \right) dX dt \nonumber \\
&= \int_0^T \int_\Omega \left[ -\omega_\rho(\rho)\, \rho + \omega(\rho) \right] (\nabla_x \cdot \delta x) \, dxdt \nonumber \\
&= \int_0^T \int_\Omega \nabla_x \left[ \omega_\rho(\rho)\, \rho - \omega(\rho) \right] \cdot \delta x \, dxdt 
= \left\langle \nabla_x \left[ \omega_\rho(\rho)\, \rho - \omega(\rho) \right], \delta x\right\rangle_{L^2_{x,t}}.
\end{align}
As an example, the choice $\omega(\rho)=\lambda\rho\ln\rho$ recovers the heat equation $\partial_t\rho = \lambda\Delta\rho$, where $\lambda$ is the diffusion coefficient.

\subsection{Application of EVMDOT to the Fokker--Planck equation}
The Fokker--Planck equation is written as
\begin{equation}
\label{eq:10}
\partial_t\rho+\nabla\cdot(a(x)\rho) = \frac{1}{2} \nabla\cdot (\sigma^2(x)\nabla\rho),
\end{equation}
where $a(x)\in\mathbb{R}^n$ is the drift term and $\sigma(x)\in\mathbb{R}$ is the noise intensity. In accordance with the fluctuation--dissipation theorem, we constrain the convection coefficient to $a(x)=-\frac{1}{2}\sigma^2(x)\nabla\psi(x)$, where $\psi(x)\in\mathbb{R}$ is the potential function. We then construct a variational model in which the free energy is
\begin{equation}
\mathcal{F}(\rho,\psi)=\int _\Omega \rho(x,t) \ln \rho(x,t)+\psi(x)\rho(x,t) dx,
\end{equation}
where $\phi(x,t) = -\ln \rho(x,t)-\psi(x) $ and the energy dissipation is
\begin{equation}
\mathcal{D}(u,\rho)=\frac{1}{2}\int_\Omega \frac{\rho(x,t)}{\sigma^2(x)/2}\|u(x,t)\|^2dx.
\end{equation}
This variational model of dynamic optimal transport is recast as
\begin{align}
&\mathrm{force}_{conservative} = \nabla\rho+\rho\nabla\psi \nonumber\\
&\mathrm{force}_{dissipative} = \frac{\rho}{\sigma^2/2} u \nonumber\\
&\mathrm{force}_{conservative} + \mathrm{force}_{dissipative} = 0 \nonumber\\
&\rho_t+\nabla\cdot(\rho u) = 0  \Leftrightarrow \rho(x,t) = \rho(X,0)/\mathrm{det} F(X,t) 
\end{align}
Here, the velocity field $u(x,t)$ is obtained from Eq.~(\ref{eq:4}) as
\begin{equation}
\label{eq:5}
u(x,t) =  - \frac{\sigma^2}{2} \nabla \ln\rho - \frac{\sigma^2}{2}\nabla \psi,
\end{equation}
Substituting $u(x,t)$ into the continuity equation recovers the Fokker--Planck equation. For an isothermal closed system, combining the first and second laws of thermodynamics yields the energy dissipation law \cite{giga2017variational}
\begin{equation}
\frac{d}{dt} E^{\mathrm{total}} = - \Delta,
\end{equation}
where $E^{\mathrm{total}}$ is the sum of the kinetic energy $\mathcal{K}$ and the total Helmholtz free energy $\mathcal{F}$, and $\Delta = 2\mathcal{D}$ is the rate of energy dissipation. In terms of the present model, this law can be written as
\begin{equation}
\label{eq:6}
\frac{d}{dt}\int_\Omega[\rho \ln\rho+\psi\rho]dx = - \int_\Omega \frac{\rho}{\sigma^2/2} \|u\|^2 dx.
\end{equation}
Since the free energy $\mathcal{F}(t)$, the velocity $u(x,t)$, and the dissipation rate are all expressed in terms of the probability density function $\rho(x,t)$, it is most natural to evaluate the loss function directly from the density data $\rho(x,t)$. Specifically, the training data are observed density datasets, denoted $\rho(x,t_k)$ for $k=1,2,\ldots,K$. The free energy functional $\mathcal{F}$ at time $t_k$ for the density $\rho(x,t)$ is written as
\begin{equation}
\mathcal{F}(t_k)= \int_\Omega \rho(x,t_k)\big(\ln{\rho}(x,t_k)+\psi_{nn}(x;\theta_{\psi})\big)\,dx,
\end{equation}
where the internal energy $\psi$ is represented by a multilayer neural network $\psi_{nn}(x;\theta_{\psi})$. Following the energy-dissipation law (Eqs.~\ref{eq:5}--\ref{eq:6}), we construct the $L^2$ loss function for the $m$-th trajectory $\{\rho^{(m)}(x,t_k)\}_{k=1}^{K}$ as
\begin{align}
    \mathcal{L}^{(m)} = \bigg(&\mathcal{F}^{(m)}(t_{K})-\mathcal{F}^{(m)}(t_{1}) +\frac{\sigma^2}{2}\int_{t_1}^{t_K}
    \int_\Omega \rho(x,t)\|\nabla\ln\rho(x,t)+\nabla\psi_{nn}(x;\theta_{\psi})\|_2^2\,dx\,dt\bigg)^2.
\end{align}
For simplicity, the noise intensity $\sigma$ is held fixed. In practice, both the spatial integrals over $\Omega\subset\mathbb{R}^{d}$ and the temporal summation are evaluated using the composite trapezoidal rule. Let $\Omega=\prod_{j=1}^{d}[a_j,b_j]$ be discretized by a tensor-product grid with $N_j+1$ nodes along the $j$-th coordinate, spacing $h_j=(b_j-a_j)/N_j$, and nodes $x_{j,i_j}=a_j+i_j h_j$ for $i_j=0,1,\ldots,N_j$. A multi-index $h=(i_1,\ldots,i_d)$ identifies the grid point $x_{\mathbf{i}}=(x_{1,i_1},\ldots,x_{d,i_d})$. The one-dimensional trapezoidal weight along direction $j$ is
\begin{equation}
\label{eq:9}
w_{j,i_j}=
\begin{cases}
\tfrac{1}{2}h_j, & i_j=0 \text{ or } i_j=N_j,\\
h_j, & 1\le i_j\le N_j-1,
\end{cases}
\end{equation}
and the $d$-dimensional composite weight is the tensor product $W_{h}=\prod_{j=1}^{d} w_{j,i_j}$. For any integrand $g:\Omega\to\mathbb{R}$, the spatial integral is approximated by
\begin{equation}
\int_{\Omega} g(x)\,dx \;\approx\; 
\sum_{h} W_{h}\, g(x_{h}).
\end{equation}
Applying this rule to the free energy and the dissipation yields the discrete formulas for the $m$-th trajectory:
\begin{align}
\label{eq:7}
\mathcal{F}^{(m)}(t_k) \;&=\; 
\sum_{h} W_h\,
\rho^{(m)}(x_h,t_k)\Big(\ln{\rho}^{(m)}(x_h,t_k)
+\psi_{nn}(x_h;\theta_{\psi})\Big), \nonumber \\
\mathcal{D}^{(m)}(t_k) \;&=\;\sum_{h} W_h\,\rho^{(m)}(x_{h},t_k)\,\Big\|\nabla\ln{\rho}^{(m)}(x_{h},t_k)+\nabla\psi_{nn}(x_h;\theta_\psi)\Big\|_2^2
\end{align}
For the time grid $\{t_k\}_{k=1}^{K}$ with non-uniform steps $\Delta t_k=t_{k+1}-t_k$, the composite trapezoidal rule applied in time to the dissipation integral $\int_{t_1}^{t_K}\tfrac{\sigma^2}{2}\mathcal{D}(t)\,dt$ gives
\begin{equation}
\sum_{k=1}^{K-1}\frac{\sigma^2}{2}\cdot\frac{\Delta t_k}{2}
\big(\mathcal{D}^{(m)}(t_k)+\mathcal{D}^{(m)}(t_{k+1})\big).
\end{equation}
where $\Delta t_k = t_{k+1}-t_k$. Combining the spatial and temporal trapezoidal rules, the discrete loss for the $m$-th trajectory becomes
\begin{align}
\label{eq:8}
    \mathcal{L}^{(m)} 
    =\Bigg(&\;\mathcal{F}^{(m)}(t_{K})-\mathcal{F}^{(m)}(t_{1}) +\sum_{k=1}^{K-1}\frac{\sigma^2}{4}\,\Delta t_k
    \Big(\mathcal{D}^{(m)}(t_k)+\mathcal{D}^{(m)}(t_{k+1})\Big)\Bigg)^2.
\end{align}
The total loss is the sum of all trajectory losses, and the neural network $\psi(x;\theta_\psi)$ is trained by minimizing
\begin{equation}
\min_{\theta_\psi}\mathcal{L} = \min_{\theta_\psi} \sum_{m=1}^{M} \lambda^{(m)}\,\mathcal{L}^{(m)}. 
\end{equation}
where $\lambda^{(m)}$ is the weight of the $m$-th trajectory.
\section{Results}
\subsection{Reconstruction of the Energy Landscape: Balancing Data Quality and Data Quantity}
\begin{figure}[ht]
\centering
\includegraphics[width=0.9\linewidth]{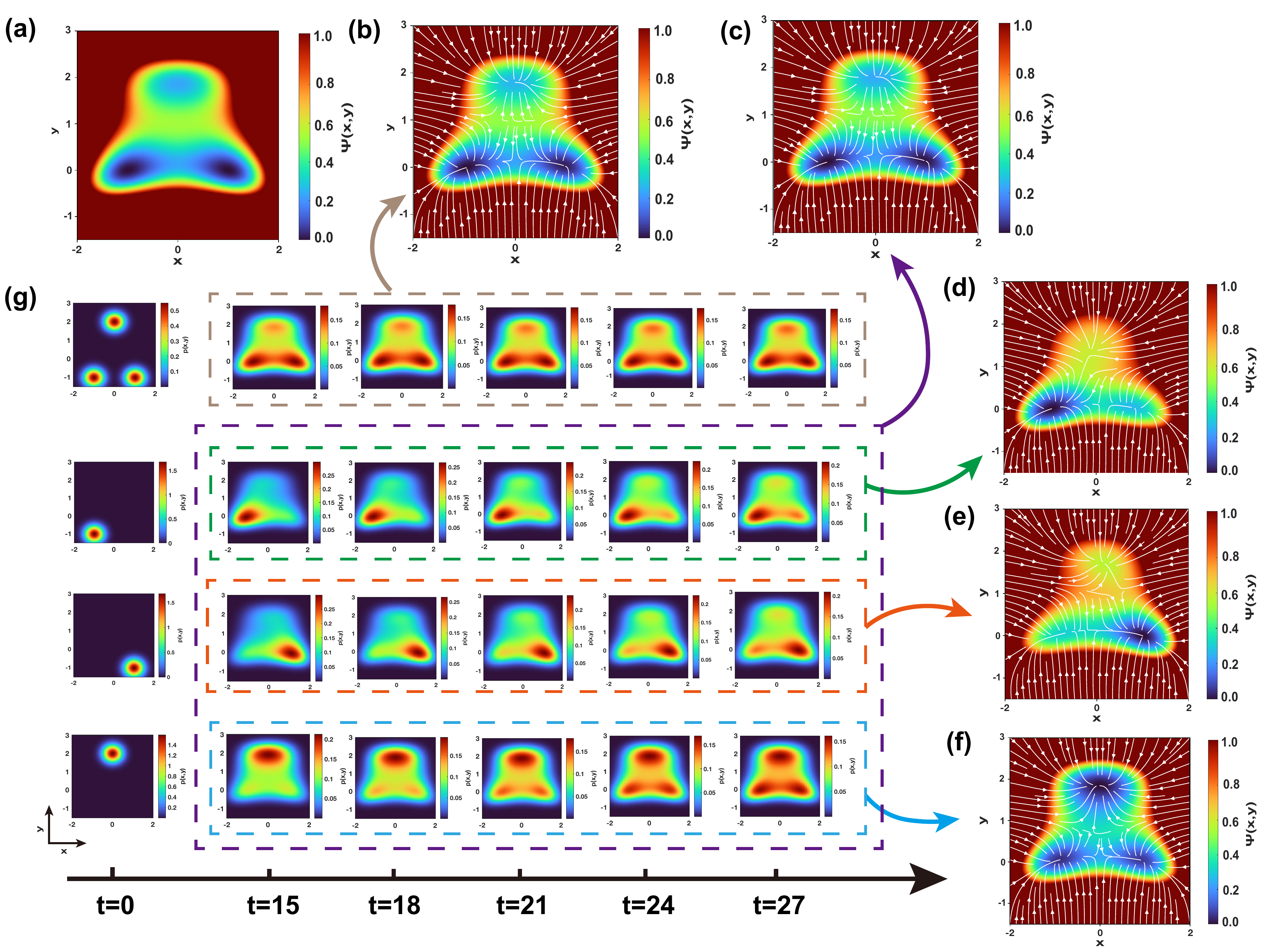}
\caption{\scriptsize\textbf{Energy landscape reconstruction.}
\textbf{(a)} Energy landscape of the benchmark model.
\textbf{(b)} Reconstructed energy landscape obtained from a trajectory initialized over all three wells.
\textbf{(c)} Reconstructed energy landscape obtained from the three trajectories in (d)--(f) combined.
\textbf{(d)} Reconstructed energy landscape obtained from a trajectory whose initial distribution is peaked near $x_A$.
\textbf{(e)} Reconstructed energy landscape obtained from a trajectory whose initial distribution is peaked near $x_B$.
\textbf{(f)} Reconstructed energy landscape obtained from a trajectory whose initial distribution is peaked near $x_C$.
The white arrows denote the vector field computed from $\psi_{nn}(x;\theta_\psi)$ as $a(x)=-\frac{1}{2}\sigma^2(x)\nabla\psi(x)$, where $\psi(x)\in\mathbb{R}$, and the colorbar indicates the value of the potential $\Psi$. 
\textbf{(g)} Simulation data generated from diverse initial conditions using the Fokker--Planck equation (Eq.~\eqref{eq:10}) together with the benchmark model (Eq.~\eqref{eq:11}); the colorbar indicates the probability density $\rho$. The value of $\sigma^2$ is set to 0.2.}
\label{fig:1}
\end{figure}

To assess the applicability of the EVMDOT, a benchmark model (Fig.~\ref{fig:1}a) was constructed as
\begin{equation}
\label{eq:11}
\psi(x,y) = (\frac{x^2}{2}+\frac{y-1}{2})^2+\frac{((y-1)^2-1)^2}{2}
\end{equation}
With a spatially uniform initial distribution, the early-phase data reveal the three wells of the energy landscape (Fig.~\ref{fig:S3}), located at
$x_A = (0,2)^\top$, $x_B = (-1,0)^\top$, and $x_C = (1,0)^\top$
. For the simulation data (See details on A.5 and Algorithm \ref{alg:fp_solver}), we restrict our attention to
the late-phase time points ($t \in [15,27]$), at which the probability
density function has approached its stationary distribution
(Fig.~\ref{fig:1}g). A no-flux boundary condition is imposed
($\nabla\rho \cdot \mathbf{n}|_{\partial\Omega} = 0$). When the initial
condition is taken as a Gaussian mixture supported on the three wells, the
energy landscape reconstructed by the EVMDOT agrees closely with the benchmark (Fig.~\ref{fig:1}b), demonstrating the reliability of the EVMDOT in the regime of sufficiently high data quality. The white arrows display the drift field $a(x)$ rather than the transport velocity $u(x,t)$, because $a(x)$ is the autonomous, time-independent vector field that characterizes the intrinsic dynamics of the system. Since the diffusive term $-\frac{\sigma^2}{2}\nabla\ln\rho(x,t)$ is already determined from the data, displaying $u(x,t)$ would add no further information.

Furthermore, when the initial condition is a single-peaked Gaussian
distribution, representing the case of insufficient data quality, the
EVMDOT fails to reconstruct the energy landscape
(Figs.~\ref{fig:1}d--f). When the peak of the initial distribution is
located near $x_A$, the deepest well of the learned energy is found at
$x_A$ (Fig.~\ref{fig:1}d); likewise, when the peak is near $x_B$ or $x_C$,
the deepest well is located at $x_B$ (Fig.~\ref{fig:1}e) or $x_C$
(Fig.~\ref{fig:1}f), respectively. However, when these three trajectories
are integrated, so that the data quantity becomes sufficient, the
EVMDOT reconstructs the benchmark model almost
perfectly (Fig.~\ref{fig:1}c), revealing the generalization capability of
our method through calibration against the available data quantity.

\subsection{Reconstruction of the Waddington Landscape: Data Quantity Compensates for Limited Early-Stage Data Quality}
The Waddington landscape originated as a visual metaphor for cellular differentiation, but has since developed into a powerful mathematical and experimental framework for studying gene-regulatory dynamics, noise, and fate decisions in development \cite{feinberg2023epigenetics}, reprogramming\cite{lang2021landscape}, and cancer \cite{wang2025mathematical}. 
In Waddington's picture, cells are represented as balls rolling down branching valleys that trace out the available differentiation paths. The valley floors correspond to attractors (stable cell states, i.e. terminal fates), while the ridges between them act as barriers that stabilize these fates and define the bifurcation points separating them \cite{ferrell2012bistability}. 

The energy landscape and the Waddington landscape differ in their prior information: for the former, the benchmark energy landscape is known a priori and serves as the ground truth, whereas for the latter, the prior is the deterministic drift $\frac{d\mathbf{x}}{dt}=a(\mathbf{x})$ of the underlying dynamics. This drift is distinct from the transport velocity field $u(\mathbf{x},t)$ of dynamic optimal transport: the two are related through Eq.~\eqref{eq:5} (i.e. $u(\mathbf{x},t) = a(\mathbf{x}) - \frac{\sigma^2}{2}\nabla\ln\rho(\mathbf{x},t)$), where $a(\mathbf{x})$ drives individual trajectories while $u(\mathbf{x},t)$ advects the probability density. Consequently, the Waddington landscape is conventionally constructed from the stationary distribution $\rho_{ss}$ ($\psi (\mathbf{x})= -\ln \rho_{ss}(\mathbf{x})$), which serves as the ground truth \cite{wang2011quantifying,shi2022energy,zhao2024epr}.

\begin{figure}[ht]
\centering
\includegraphics[width=1.0\linewidth]{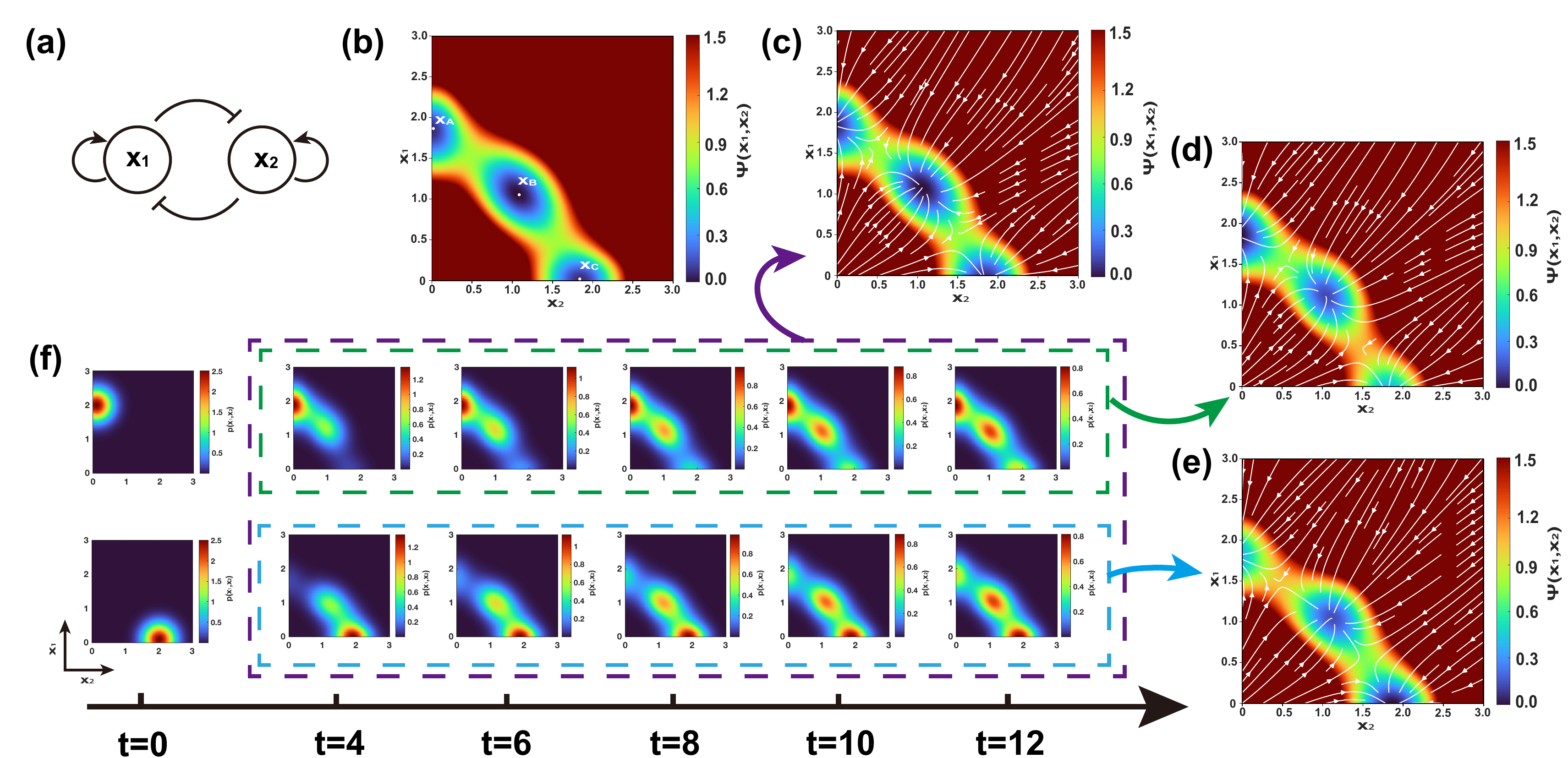}
\caption{\scriptsize\textbf{Waddington landscape reconstruction of a positive-feedback gene-regulatory motif.}
\textbf{(a)} Mutual regulation of the two genes $x_1$ and $x_2$ in the cell-fate decision system.
\textbf{(b)} Waddington landscape constructed from the stationary distribution $\rho_{ss}$ ($\psi = -\ln \rho_{ss}$).
\textbf{(c)} Reconstructed Waddington landscape obtained from the two trajectories in (d) and (e) combined.
\textbf{(d)} Reconstructed Waddington landscape obtained from a trajectory whose initial distribution is peaked near $x_A$.
\textbf{(e)} Reconstructed Waddington landscape obtained from a trajectory whose initial distribution is peaked near $x_C$.
The white arrows denote the vector field computed from $\psi_{nn}(x;\theta_\psi)$ as $a(x)=-\frac{1}{2}\sigma^2(x)\nabla\psi(x)$, where $\psi(x)\in\mathbb{R}$, and the colorbar indicates the value of the potential $\Psi$.
\textbf{(f)} Simulation data generated from diverse initial conditions using the Fokker--Planck equation (Eq.~\eqref{eq:10}) combined with the drift term (Eq.~\eqref{eq:12}); the colorbar indicates the probability density $\rho$. The value of $\sigma^2$ is set to $0.2$.}
\label{fig:2}
\end{figure}

Here, we first employ a dynamical model (Eq.~\eqref{eq:12}) of a positive-feedback gene-regulatory motif (Fig. \ref{fig:2}a)  to generate time-resolved probability density data (Fig. \ref{fig:2}f). We then reconstruct the Waddington landscape from these data using our proposed EVMDOT.
\begin{align}
\label{eq:12}
\frac{dx_1}{dt} = \frac{a_1x_1^n}{S^n+x_1^n}+\frac{b_1S^n}{S^n+x_2^n}-k_1x_1 \nonumber \\
\frac{dx_2}{dt} = \frac{a_2x_2^n}{S^n+x_2^n}+\frac{b_2S^n}{S^n+x_1^n}-k_2x_2 
\end{align}
The parameters are set to $a_1=a_2=b_1=b_2=k_1=k_2=1$, $n=4$, and $S=0.5$ \cite{wang2011quantifying,zhao2024epr}. This model features three wells on the Waddington landscape, located at $x_A$, $x_B$, and $x_C$ (Fig.~\ref{fig:2}b). For the simulation data (see Appendix~A.5 and Algorithm~\ref{alg:fp_solver}), we restrict our attention to the early-phase time points ($t \in [4,12]$), at which the probability density function has not yet reached its stationary distribution (Fig.~\ref{fig:2}f). A no-flux boundary condition is imposed ($\nabla\rho \cdot \mathbf{n}|_{\partial\Omega} = 0$). To remain consistent with the biological meaning of the Waddington landscape , which describes the differentiation of cells from one type into another \cite{ferrell2012bistability}, we set the initial condition to a Gaussian distribution concentrated in a single well. This choice mirrors the developmental scenario in which a population of cells initially occupies a single progenitor state and subsequently differentiates into distinct terminal fates as it evolves across the landscape \cite{wang2011quantifying}. 

When trained on the data whose initial distribution is peaked near $x_A$, our method accurately identifies the depths of the wells at $x_A$ and $x_B$, but fails to resolve the depth of the well at $x_C$ (Fig. \ref{fig:2}d). 
Conversely, when trained on the data whose initial distribution is peaked near $x_C$, our method accurately identifies the depths of the wells at $x_C$ and $x_B$, but fails to resolve the depth of the well at $x_A$ (Fig. \ref{fig:2}e). This behavior is consistent with the data themselves: trajectories generated from different initial conditions differ in their density distribution over the three wells during the early-stage dynamics (Fig. \ref{fig:2}f). However, when these two trajectories are integrated, our method reconstructs the Waddington landscape perfectly (Fig. \ref{fig:2}c), in full agreement with the ground truth (Fig. \ref{fig:2}b). This indicates that our EVMDOT strikes a balance between data quality and data quantity: the insufficient data quality caused by limited temporal sampling can be compensated for by sufficient data quantity, thereby enabling faithful reconstruction of the Waddington landscape.

\begin{align}
\label{eq:13}
\frac{dx_1}{dt} &= \frac{\alpha_1 x_1^2 + \beta}{1 + \alpha_1 x_1^2 + \gamma_2 x_2^2 + \gamma_3 x_3^2 + \beta} - \delta_1 x_1, \nonumber\\
\frac{dx_2}{dt} &= \frac{\alpha_2 x_2^2 + \beta}{1 + \gamma_1 x_1^2 + \alpha_2 x_2^2 + \gamma_3 x_3^2 + \beta} - \delta_2 x_2, \nonumber\\
\frac{dx_3}{dt} &= \frac{\alpha_3 x_3^2}{1 + \alpha_3 x_3^2} - \delta_3 x_3.
\end{align}

\begin{figure}[ht]
\centering
\includegraphics[width=1.0\linewidth]{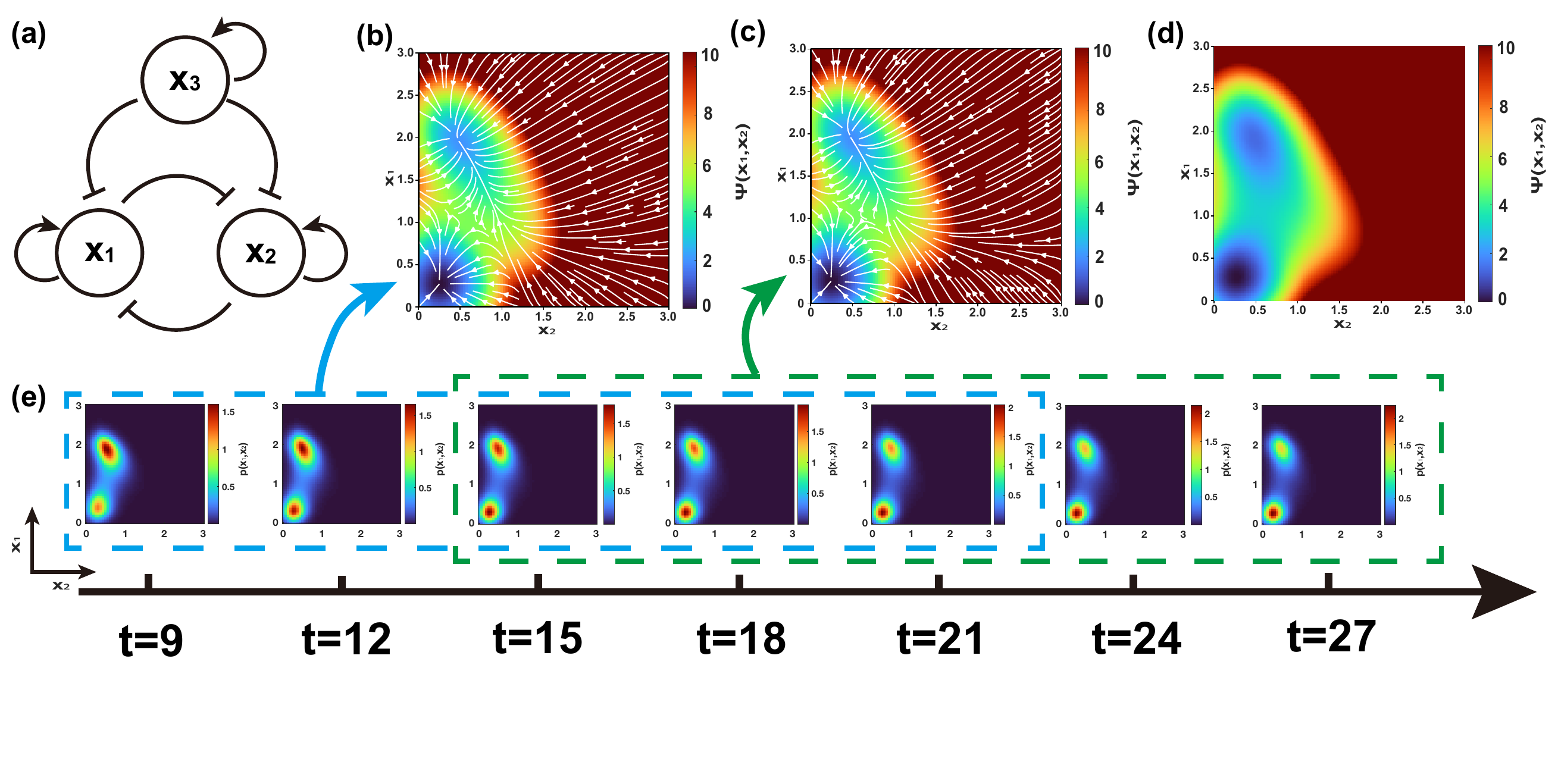}
\caption{\scriptsize\textbf{Waddington landscape reconstruction of a three-gene regulatory motif.}
\textbf{(a)} Schematic diagram of the three-gene regulatory model.
\textbf{(b)} Reconstructed Waddington landscape obtained from a trajectory in the early phase.
\textbf{(c)} Reconstructed Waddington landscape obtained from a trajectory in the late phase.
\textbf{(d)} Waddington landscape constructed from the stationary distribution $\rho_{ss}$ ($\psi = -\ln \rho_{ss}$). The white arrows denote the vector field computed from $\psi_{nn}(x;\theta_\psi)$ as $a(x)=-\frac{1}{2}\sigma^2(x)\nabla\psi(x)$, where $\psi(x)\in\mathbb{R}$, and the colorbar indicates the value of the potential $\Psi$.
\textbf{(e)} Simulation data generated from diverse initial conditions using the Fokker--Planck equation (Eq.~\eqref{eq:10}) combined with the benchmark model (Eq.~\eqref{eq:13}); the colorbar indicates the probability density $\rho$. The value of $\sigma^2$ is set to $0.02$.}
\label{fig:3}
\end{figure}

Furthermore, we apply the EVMDOT to reconstruct the Waddington landscape of a three-gene motif (Fig. \ref{fig:3}a). The dynamical model of this three-gene motif \cite{sun2026variational,sha2024reconstructing} and the corresponding parameter settings are detailed in Eq.~\eqref{eq:13} and Table~\ref{tab:1} \cite{sun2026variational}, respectively. To facilitate visualization in this case, we define the projection potential \cite{lv2015energy} $\Psi(x_1,x_2) = \min_{x_3}\Psi(x_1,x_2,x_3)$ , obtained by minimizing the full three-dimensional potential over $x_3$ .

Owing to the additional inhibitory regulation of $x_1$ and $x_2$ by gene $x_3$, the Waddington landscape degenerates into two wells, and the symmetry between the wells is also broken. By combining the stationary distribution with our definition of the reduced potential, the ground truth can still be constructed (Fig.~\ref{fig:3}d). For the simulation data (See details on A.5 and Algorithm \ref{alg:fp_solver}), we generate $64$ distinct trajectories from $64$ different initial distributions, each a Gaussian centered at a random point with the same fixed variance. The time points are selected over $t \in [9,27]$, and a representative trajectory is shown in Fig.~\ref{fig:3}e. A no-flux boundary condition is imposed ($\nabla\rho \cdot \mathbf{n}|_{\partial\Omega} = 0$). We further partition the generated trajectories into two phases: $t \in [9,21]$ corresponds to the early phase, and $t \in [15,27]$ to the late phase.

The results show that the Waddington landscape reconstructed from the early-phase data differs only marginally from that reconstructed from the late-phase data, with both being consistent with the ground truth (Figs.~\ref{fig:3}b-c). This further corroborates the mechanism identified above: when the data quantity is sufficient, our EVMDOT strikes a balance between data quality and data quantity, rendering it robust to variations in data quality that arise from the choice of the temporal window.

\subsection{Inference of the potential landscape of A$\beta$ and tau in Alzheimer's disease}
Alzheimer's disease (AD) is the leading cause of dementia worldwide and poses a major and growing public health burden. Its hallmark pathology involves the accumulation of amyloid-$\beta$ (A$\beta$) plaques and the progressive spread of misfolded tau aggregates, which are closely associated with synaptic dysfunction, neurodegeneration, and cognitive decline. A central open question is how the coupled dynamics of A$\beta$ and tau generate the observed spatiotemporal progression of disease across heterogeneous populations \cite{hao2016mathematical}.
Recent work in computational neuroscience and mathematical biology suggests that AD progression can be described using dynamical systems and reaction--diffusion-type models capturing the interaction and propagation of pathological proteins \cite{hao2022optimal}. More recent studies further develop mechanistic and data-driven frameworks for A$\beta$--tau coupling using longitudinal neuroimaging data from ADNI, highlighting the role of their interaction in shaping disease trajectories and revealing latent dynamical structures underlying biomarker evolution \cite{wang2026learning,zheng2022data}.

Building on this growing body of work, 
We apply the EVMDOT to the Alzheimer's Disease Neuroimaging Initiative (ADNI) dataset to infer the potential landscape of A$\beta$ and tau. Unlike the previous cases, in which the probability density data are generated from the Fokker--Planck equation, the ADNI dataset requires preprocessing to be converted into time-series probability density data, which then serve as the input to the EVMDOT for inferring the potential. 

\begin{figure}[ht]
\centering
\includegraphics[width=1.0\linewidth]{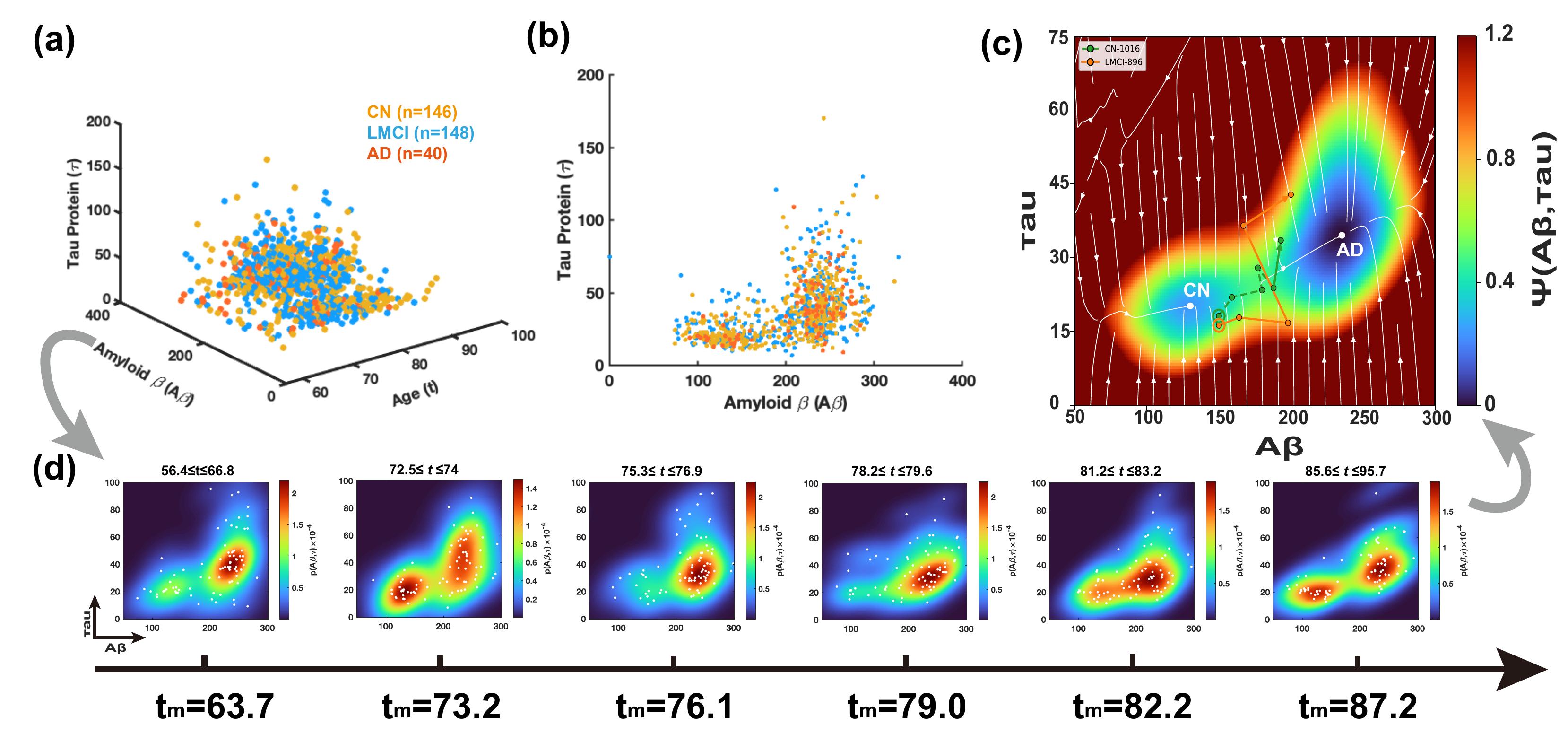}
\caption{\scriptsize\textbf{Inference of the potential landscape of A$\beta$ and tau.}
\textbf{(a)} Scatter plot of the spatio-temporal distribution of A$\beta$ and tau. Orange, blue, and red dots represent participants in the CN, LMCI, and AD stages, respectively.
\textbf{(b)} Scatter projection on the spatial distribution of A$\beta$ and tau.
\textbf{(c)} Inferred potential landscape obtained using the EVMDOT. In the legend, each number denotes the participant index. The trajectories of the two participant groups are distinguished by color: green for the CN group and orange for the LMCI group. Open circles mark the starting points, and the arrowheads at the trajectory ends indicate the direction of time. The white arrows denote the vector field computed from $\psi_{nn}(x;\theta_\psi)$ as $a(x)=-\frac{1}{2}\sigma^2(x)\nabla\psi(x)$, where $\psi(x)\in\mathbb{R}$, and the colorbar indicates the value of the potential $\Psi$. The value of $\sigma^2$ is set to $1.0$.
\textbf{(d)} Time-series probability density distributions obtained by applying two-dimensional kernel density estimation to the ADNI dataset. Here, $t_m$ denotes the median of the $m$-th time interval, and the colorbar indicates the probability density $\rho$.
}
\label{fig:4}
\end{figure}

The ADNI dataset comprises $146$ cognitively normal (CN) participants, $148$ participants with late mild cognitive impairment (LMCI), and $40$ participants with Alzheimer's disease (AD)  \cite{hao2022optimal}, each of whom was followed longitudinally with measurements recorded at multiple time points (Fig.~\ref{fig:4}a). Because each participant has longitudinal data, it is difficult to find data from multiple participants sharing exactly the same time point. We therefore divide the time axis into $12$ intervals, each containing $N_m = 80$ data points pooled across participants. This preprocessing reduces the computational error caused by overly long time intervals, while ensuring that sufficient points are available in the A$\beta$--tau space to estimate the probability density $\hat{\rho}(\mathbf{x}, t_m)$ (Fig.~\ref{fig:4}b) via two-dimensional kernel density estimation (Eq.~\eqref{eq:14}).
 \begin{equation}
\label{eq:14}
\hat{\rho}(\mathbf{x}, t_m) = \frac{1}{N_m}\sum_{i=1}^{N_m} 
\frac{1}{2\pi h_1 h_2} 
\exp\!\left(
-\frac{(x_1 - x_1^{(i)})^2}{2 h_1^2}
-\frac{(x_2 - x_2^{(i)})^2}{2 h_2^2}
\right).
\end{equation}
 The bandwidth is selected using Scott's rule, the default method in the SciPy implementation of Gaussian kernel density estimation (\texttt{scipy.stats.gaussian\_kde}). For each interval, the median of the time values is used to approximate the representative time point of the corresponding data. For clarity, Fig.~\ref{fig:4}d shows the probability density distributions for only $6$ of the $12$ time intervals. 
 
To incorporate the density structure of the data as a prior that guides the optimization direction of the EVMDOT, we augment the original loss function with a regularization term (Eq.~\eqref{eq:15}):
\begin{align}
\label{eq:15}
\min_{\theta_\psi}\mathcal{L} 
= \min_{\theta_\psi} \left[ \mathcal{L}^{(1)} 
+ \lambda \sum_{m=1}^{N_T} \int_{\Omega} \rho(\mathbf{x}, t_m)\big(\Delta_m - \mathbb{E}[\Delta_m]\big)^2 \, d\mathbf{x} \right],\nonumber \\
\Delta_m = \psi_{nn}(\mathbf{x};\theta_\psi) - \ln \rho(\mathbf{x}, t_m),
\qquad
\mathbb{E}[\Delta_m] = \int_{\Omega} \rho(\mathbf{x}, t_m)\, \Delta_m \, d\mathbf{x}.
\end{align}
The regularization term penalizes the density-weighted variance of $\Delta_m = \psi_{nn}(\mathbf{x};\theta_\psi) - \ln\rho(\mathbf{x}, t_m)$ about its mean $\mathbb{E}[\Delta_m]$ over each time interval. Minimizing this term encourages $\psi_{nn}(\mathbf{x};\theta_\psi) - \ln\rho(\mathbf{x}, t_m)$ to remain approximately constant across the high-density region, thereby softly enforcing the relation $\psi_{nn}(\mathbf{x};\theta_\psi) \approx -\ln\rho(\mathbf{x}, t_m) + \mathrm{const}$. This is consistent with the stationary relation $\psi = -\ln\rho_{ss}$, and the weight $\lambda=0.1$ controls the strength of the regularization.

The potential landscape inferred by the EVMDOT exhibits two wells (Fig.~\ref{fig:4}c), in agreement with the data distribution that serves as its input (Fig.~\ref{fig:4}b). The well in which both A$\beta$ and tau are at low levels can be interpreted as the CN stage, whereas the well in which both are at high levels corresponds to the AD stage. To validate this interpretation independently of the density data used for inference, we examined individual participants' longitudinal trajectories on the inferred landscape, which encode the temporal direction of progression rather than the static density. We found two trajectories transitioning from the CN well to the AD well, both passing through the saddle point on the ridge between the two wells (Fig.~\ref{fig:4}c). This matches the expected escape pathway over the lowest barrier between the two attractors \cite{weinan2002string,lv2014constructing,yin2021searching}, supporting the validity of the reconstructed landscape. Importantly, the CN well is shallower and is separated from the AD well by a low potential barrier, indicating that the CN state can readily escape its well under an external driving force and transition into the AD state (Fig.~\ref{fig:4}c). This confirms that the EVMDOT, applied to the ADNI data, captures the biological reality of disease-state progression in participants, together with the key transition pathway between states.
\section{Discussions}
Dynamic optimal transport reformulates optimal transport as the evolution of densities and velocities over time. Following the Benamou--Brenier fluid-mechanics formulation\cite{benamou2000computational,benamou2002monge}, the static transport is replaced by a time-dependent density $\rho(x,t)$ and velocity field $u(x,t)$ satisfying the continuity equation $\partial_t\rho + \nabla\cdot(\rho u) = 0$; the optimal coupling is then obtained by minimizing the kinetic energy. In this paper, we introduce the flow map and mechanics to reformulate dynamic optimal transport within an energetic variational framework. For the Fokker--Planck equation, we develop the EVMDOT to reconstruct the energy landscape and the Waddington landscape from time-series density data $\rho(x,t)$. Furthermore, using the ADNI dataset, we apply the EVMDOT to infer the potential landscape of A$\beta$ and tau, thereby identifying the disease-state progression of Alzheimer's disease.

The flow map provides a natural starting point. From a mechanical perspective, it is a coordinate transformation between the Eulerian and Lagrangian descriptions. Under the Lagrangian description, the density is given by $\rho(x,t) = \rho(X,0)/\det F(X,t)$ through the deformation gradient $F(X,t)$, so that the continuity-equation constraint imposed in the dynamic optimal transport formulation is satisfied automatically. The originally constrained optimization is thereby recast as an unconstrained one. In other words, transforming the frame of reference in which the observed fluid is described is equivalent to enforcing the continuity equation satisfied by the fluid. 

In mechanics, the optimization objective of dynamic optimal transport is equivalent to the force balance jointly characterized by the least action principle (LAP) and the maximum dissipation principle (MDP), and the Lagrange multiplier introduced in the Benamou--Brenier formulation is precisely the action functional. Specifically, the action functional $\mathcal{A}$ gives rise to the conservative force through the LAP, while the dissipation functional $\mathcal{D}(x)$ gives rise to the dissipative force through the MDP; the balance between these two forces yields the relation $u(x,t) = \nabla\phi(x,t)$. Combining this with the flow map, we obtain the Hamilton--Jacobi equation derived from the Benamou--Brenier formulation. In this way, dynamic optimal transport is seamlessly recast as an energetic variational model.

The EVMDOT transforms the constrained optimization problem into an unconstrained one, in which the velocity field $u(x,t)$ is obtained by balancing the conservative and dissipative forces, obtained by taking variations of the free energy $\mathcal{F}(x)$ and the dissipation $\mathcal{D}(x)$. Moreover, this framework offers excellent extensibility: by specifying different forms of the free energy, one can derive the corresponding velocity field $u(x,t)$ \cite{giga2017variational}, although in this study we consider only the free-energy form associated with the Fokker--Planck equation \cite{lu2024learning}. From a complementary viewpoint, the same variational structure of the Fokker--Planck equation can be recovered through the theory of gradient flows \cite{hu2024energetic} in the Wasserstein space, where the dynamics are characterized as the steepest descent of the free energy. In particular, the JKO scheme \cite{jordan1998variational,fu2023high} realizes this gradient flow as a sequence of variational problems involving the Wasserstein distance, providing a time-discrete counterpart to the continuous energetic variational formulation adopted here.

Many existing methods for reconstructing the energy landscape rely on sampling data that have approached the stationary distribution \cite{shi2022energy,zhao2024epr,lu2024learning}. In this regime, the data are of high quality in the sense that they conform to the theoretical relation $\psi(x) = -\ln \rho_{ss}(x)$, where $\rho_{ss}$ denotes the stationary density; the energy landscape can then be recovered simply by taking the negative logarithm of the empirical density, yielding a reconstruction consistent with the ground truth. However, this requirement is often restrictive in practice, as reaching the stationary distribution may demand prohibitively long sampling times, and the observed data are frequently confined to transient, non-stationary stages of the dynamics. This limitation motivates us to examine how the quality and the quantity of the sampling data jointly affect the reconstruction. In reconstructing the energy landscape with the EVMDOT, we account for both data quality and data quantity in the design of the data. In our experiments, high data quality corresponds to the optimal sampling trajectory, i.e. the data generated from the optimal initial condition, whereas high data quantity corresponds to multiple distinct sampling trajectories, each of insufficient quality on its own (Fig.~\ref{fig:1}g). In practical applications, data quantity can compensate for limited data quality, and is often the more decisive factor. The EVMDOT successfully reconstructs the energy landscape from high-quality data (Fig.~\ref{fig:1}b), confirming the validity of the method itself; meanwhile, high-quantity data more clearly demonstrate the ability of the EVMDOT to extract the informative features from the data and to reconstruct the energy landscape accurately (Fig.~\ref{fig:1}c).

In the numerical experiments reconstructing the Waddington landscape with the EVMDOT, we deliberately generated data of low quality induced by the temporal sampling (Figs~\ref{fig:2}f and~\ref{fig:3}e). Even so, increasing the data quantity was sufficient for the EVMDOT to recover the Waddington landscape accurately from such degraded data. This observation highlights a complementary relationship between data quality and data quantity within the EVMDOT: a shortfall in one can be offset by an abundance of the other, so that the choice of the temporal window has little impact on the reconstruction. However, in our studies, the Waddington landscapes constructed by the EVMDOT are all multi-well cases that can be approximated as gradient systems; the EVMDOT cannot construct a Waddington landscape with circulating flux. This is because a gradient system corresponds to detailed balance \cite{zhou2016construction}, in which the dynamics are driven solely by the gradient of the potential. For non-equilibrium systems, by contrast, the probability flux generally decomposes into a potential-gradient component and a non-vanishing rotational (curl) component \cite{wang2010potential}. Capturing such systems would therefore require augmenting the framework with a non-gradient, rotational-flux term.

Using the ADNI data with high heterogeneity, the EVMDOT captures the biological reality of disease-state progression in participants from CN to AD, together with the key transition pathway between states (Fig.~\ref{fig:4}). One limitation, however, arises from the purely data-driven nature of this case: unlike the previous simulation data, which intrinsically obey the underlying Fokker--Planck equation, the inferred landscape exhibits a descending gradient from the high-tau toward the low-tau region, which is inconsistent with the known dynamical characteristics of tau protein \cite{jack2024revised}. Therefore, future inference of the potential landscape of A$\beta$ and tau need to incorporate biological-mechanism constraints to improve the accuracy and plausibility of the inferred results.

In summary, the proposed EVMDOT establishes a principled bridge between dynamic optimal transport and energetic variational mechanics, recasting density evolution as a force-balance problem governed by the interplay of free energy and dissipation. Beyond a specific algorithm, this perspective offers a unifying language in which landscape reconstruction, optimal transport, and non-equilibrium dynamics are treated within a single variational structure. Crucially, our results suggest that the fidelity of a reconstructed landscape is governed less by the perfection of any single observation than by the collective coverage of the data, pointing to a data-economical principle for inferring dynamics from sparse, heterogeneous measurements. By extending this framework from controlled simulations to the clinical heterogeneity of the ADNI cohort, we demonstrate that energetic variational modeling can transform time-resolved snapshots of a complex biological system into an interpretable potential landscape, opening a quantitative route toward understanding disease progression as motion on an inferred energetic terrain.
\section{Acknowledgment}
\subsection{Funding}
S.W. and W.H. were supported by National Institute of General Medical Sciences through grant 1R35GM146894 and the Huck Chair in AI Mathematical Modeling from Penn State University's Huck Institutes of the Life Sciences. C.L. was partially supported by National Science Foundation through two grants DMS-2410742 and DMS-2118181.
\subsection{Author contributions}
W.H. and C.L. designed research; S.W. performed research; C.L. and W.H. contributed new reagents/analytic tools; S.W. analyzed data; and S.W., C.L. and W.H. wrote the paper.
\subsection{Competing interests}
There are no competing interests to declare.
\subsection{Data and code availability:} 
ADNI datasets used in this article were obtained from the ADNI database  (http://adni.loni.usc.edu/). 
Source codes and data have been deposited on the GitHub repository (https://github.com/WilliamMoriaty/EVMDOT) .

\bibliography{sample} 

@article{benamou2002monge,
  title={The Monge--Kantorovitch mass transfer and its computational fluid mechanics formulation},
  author={Benamou, J-D and Brenier, Yann and Guittet, Kevin},
  journal={International Journal for Numerical methods in fluids},
  volume={40},
  number={1-2},
  pages={21--30},
  year={2002},
  publisher={Wiley Online Library}
}

@article{benamou2000computational,
  title={A computational fluid mechanics solution to the Monge-Kantorovich mass transfer problem},
  author={Benamou, Jean-David and Brenier, Yann},
  journal={Numerische Mathematik},
  volume={84},
  number={3},
  pages={375--393},
  year={2000},
  publisher={Springer-Verlag Berlin/Heidelberg}
}

@article{hu2024energetic,
  title={Energetic variational neural network discretizations of gradient flows},
  author={Hu, Ziqing and Liu, Chun and Wang, Yiwei and Xu, Zhiliang},
  journal={SIAM Journal on Scientific Computing},
  volume={46},
  number={4},
  pages={A2528--A2556},
  year={2024},
  publisher={SIAM}
}

@article{jordan1998variational,
  title={The variational formulation of the Fokker--Planck equation},
  author={Jordan, Richard and Kinderlehrer, David and Otto, Felix},
  journal={SIAM journal on mathematical analysis},
  volume={29},
  number={1},
  pages={1--17},
  year={1998},
  publisher={SIAM}
}

@article{giga2017variational,
  title={Variational modeling and complex fluids},
  author={Giga, Mi-Ho and Kirshtein, Arkadz and Liu, Chun},
  journal={Handbook of mathematical analysis in mechanics of viscous fluids},
  pages={1--41},
  year={2017},
  publisher={Springer International Publishing, Cham}
}

@article{lu2024learning,
  title={Learning Generalized Diffusions using an Energetic Variational Approach},
  author={Lu, Yubin and Li, Xiaofan and Liu, Chun and Tang, Qi and Wang, Yiwei},
  journal={arXiv preprint arXiv:2412.04480},
  year={2024}
}

@book{su1999discontinuous,
  title={Discontinuous solutions in L (infinity) of Hamilton-Jacobi equations},
  author={Su, Bo},
  year={1999},
  publisher={Northwestern University}
}

@article{wang2011quantifying,
  title={Quantifying the Waddington landscape and biological paths for development and differentiation},
  author={Wang, Jin and Zhang, Kun and Xu, Li and Wang, Erkang},
  journal={Proceedings of the National Academy of Sciences},
  volume={108},
  number={20},
  pages={8257--8262},
  year={2011},
  publisher={National Academy of Sciences}
}

@article{shi2022energy,
  title={Energy landscape decomposition for cell differentiation with proliferation effect},
  author={Shi, Jifan and Aihara, Kazuyuki and Li, Tiejun and Chen, Luonan},
  journal={National Science Review},
  volume={9},
  number={8},
  pages={nwac116},
  year={2022},
  publisher={Oxford University Press}
}

@article{zhao2024epr,
  title={EPR-Net: constructing a non-equilibrium potential landscape via a variational force projection formulation},
  author={Zhao, Yue and Zhang, Wei and Li, Tiejun},
  journal={National Science Review},
  volume={11},
  number={7},
  pages={nwae052},
  year={2024},
  publisher={Oxford University Press}
}

@article{sha2024reconstructing,
  title={Reconstructing growth and dynamic trajectories from single-cell transcriptomics data},
  author={Sha, Yutong and Qiu, Yuchi and Zhou, Peijie and Nie, Qing},
  journal={Nature Machine Intelligence},
  volume={6},
  number={1},
  pages={25--39},
  year={2024},
  publisher={Nature Publishing Group UK London}
}

@article{lv2015energy,
  title={Energy landscape reveals that the budding yeast cell cycle is a robust and adaptive multi-stage process},
  author={Lv, Cheng and Li, Xiaoguang and Li, Fangting and Li, Tiejun},
  journal={PLoS computational biology},
  volume={11},
  number={3},
  pages={e1004156},
  year={2015},
  publisher={Public Library of Science San Francisco, CA USA}
}

@article{sun2026variational,
  title={Variational regularized unbalanced optimal transport: Single network, least action},
  author={Sun, Yuhao and Zhang, Zhenyi and Wang, Zihan and Li, Tiejun and Zhou, Peijie},
  journal={Advances in Neural Information Processing Systems},
  volume={38},
  pages={16963--17005},
  year={2026}
}

@article{hao2022optimal,
  title={Optimal anti-amyloid-beta therapy for Alzheimer’s disease via a personalized mathematical model},
  author={Hao, Wenrui and Lenhart, Suzanne and Petrella, Jeffrey R},
  journal={PLoS computational biology},
  volume={18},
  number={9},
  pages={e1010481},
  year={2022},
  publisher={Public Library of Science San Francisco, CA USA}
}

@article{zhou2016construction,
  title={Construction of the landscape for multi-stable systems: Potential landscape, quasi-potential, A-type integral and beyond},
  author={Zhou, Peijie and Li, Tiejun},
  journal={The Journal of chemical physics},
  volume={144},
  number={9},
  year={2016},
  publisher={AIP Publishing}
}

@article{wang2010potential,
  title={Potential and flux landscapes quantify the stability and robustness of budding yeast cell cycle network},
  author={Wang, Jin and Li, Chunhe and Wang, Erkang},
  journal={Proceedings of the National Academy of Sciences},
  volume={107},
  number={18},
  pages={8195--8200},
  year={2010},
  publisher={National Academy of Sciences}
}

@article{jack2024revised,
  title={Revised criteria for diagnosis and staging of Alzheimer's disease: Alzheimer's Association Workgroup},
  author={Jack Jr, Clifford R and Andrews, J Scott and Beach, Thomas G and Buracchio, Teresa and Dunn, Billy and Graf, Ana and Hansson, Oskar and Ho, Carole and Jagust, William and McDade, Eric and others},
  journal={Alzheimer's \& Dementia},
  volume={20},
  number={8},
  pages={5143--5169},
  year={2024},
  publisher={Wiley Online Library}
}

@article{van2023thermodynamic,
  title={Thermodynamic unification of optimal transport: Thermodynamic uncertainty relation, minimum dissipation, and thermodynamic speed limits},
  author={Van Vu, Tan and Saito, Keiji},
  journal={Physical Review X},
  volume={13},
  number={1},
  pages={011013},
  year={2023},
  publisher={APS}
}

@article{papadakis2014optimal,
  title={Optimal transport with proximal splitting},
  author={Papadakis, Nicolas and Peyr{\'e}, Gabriel and Oudet, Edouard},
  journal={SIAM Journal on Imaging Sciences},
  volume={7},
  number={1},
  pages={212--238},
  year={2014},
  publisher={SIAM}
}

@inproceedings{tong2020trajectorynet,
  title={Trajectorynet: A dynamic optimal transport network for modeling cellular dynamics},
  author={Tong, Alexander and Huang, Jessie and Wolf, Guy and Van Dijk, David and Krishnaswamy, Smita},
  booktitle={International conference on machine learning},
  pages={9526--9536},
  year={2020},
  organization={PMLR}
}

@article{burger2023dynamic,
  title={Dynamic optimal transport on networks},
  author={Burger, Martin and Humpert, Ina and Pietschmann, Jan-Frederik},
  journal={ESAIM: Control, Optimisation and Calculus of Variations},
  volume={29},
  pages={54},
  year={2023},
  publisher={EDP Sciences}
}

@article{lavenant2018dynamical,
  title={Dynamical optimal transport on discrete surfaces},
  author={Lavenant, Hugo and Claici, Sebastian and Chien, Edward and Solomon, Justin},
  journal={ACM Transactions on Graphics (TOG)},
  volume={37},
  number={6},
  pages={1--16},
  year={2018},
  publisher={ACM New York, NY, USA}
}

@article{bunne2024optimal,
  title={Optimal transport for single-cell and spatial omics},
  author={Bunne, Charlotte and Schiebinger, Geoffrey and Krause, Andreas and Regev, Aviv and Cuturi, Marco},
  journal={Nature Reviews Methods Primers},
  volume={4},
  number={1},
  pages={58},
  year={2024},
  publisher={Nature Publishing Group UK London}
}

@article{ghoussoub2018optimal,
  title={Optimal transport with controlled dynamics and free end times},
  author={Ghoussoub, Nassif and Kim, Young-Heon and Palmer, Aaron Zeff},
  journal={SIAM Journal on Control and Optimization},
  volume={56},
  number={5},
  pages={3239--3259},
  year={2018},
  publisher={SIAM}
}

@article{haasler2024scalable,
  title={Scalable computation of dynamic flow problems via multimarginal graph-structured optimal transport},
  author={Haasler, Isabel and Ringh, Axel and Chen, Yongxin and Karlsson, Johan},
  journal={Mathematics of Operations Research},
  volume={49},
  number={2},
  pages={986--1011},
  year={2024},
  publisher={INFORMS}
}

@article{chen2021optimal,
  title={Optimal transport in systems and control},
  author={Chen, Yongxin and Georgiou, Tryphon T and Pavon, Michele},
  journal={Annual Review of Control, Robotics, and Autonomous Systems},
  volume={4},
  number={1},
  pages={89--113},
  year={2021},
  publisher={Annual Reviews}
}

@article{schiebinger2019optimal,
  title={Optimal-transport analysis of single-cell gene expression identifies developmental trajectories in reprogramming},
  author={Schiebinger, Geoffrey and Shu, Jian and Tabaka, Marcin and Cleary, Brian and Subramanian, Vidya and Solomon, Aryeh and Gould, Joshua and Liu, Siyan and Lin, Stacie and Berube, Peter and others},
  journal={Cell},
  volume={176},
  number={4},
  pages={928--943},
  year={2019},
  publisher={Elsevier}
}

@article{chizat2018interpolating,
  title={An interpolating distance between optimal transport and Fisher--Rao metrics},
  author={Chizat, Lenaic and Peyr{\'e}, Gabriel and Schmitzer, Bernhard and Vialard, Fran{\c{c}}ois-Xavier},
  journal={Foundations of Computational Mathematics},
  volume={18},
  number={1},
  pages={1--44},
  year={2018},
  publisher={Springer}
}

@article{chen2018neural,
  title={Neural ordinary differential equations},
  author={Chen, Ricky TQ and Rubanova, Yulia and Bettencourt, Jesse and Duvenaud, David K},
  journal={Advances in neural information processing systems},
  volume={31},
  year={2018}
}

@book{ames2014numerical,
  title={Numerical methods for partial differential equations},
  author={Ames, William F},
  year={2014},
  publisher={Academic press}
}

@article{fu2023high,
  title={High order spatial discretization for variational time implicit schemes: Wasserstein gradient flows and reaction-diffusion systems},
  author={Fu, Guosheng and Osher, Stanley and Li, Wuchen},
  journal={Journal of Computational Physics},
  volume={491},
  pages={112375},
  year={2023},
  publisher={Elsevier}
}

@article{hao2016mathematical,
  title={Mathematical model on Alzheimer’s disease},
  author={Hao, Wenrui and Friedman, Avner},
  journal={BMC systems biology},
  volume={10},
  number={1},
  pages={108},
  year={2016},
  publisher={Springer}
}

@article{wang2026learning,
  title={Learning patient-specific spatial biomarker dynamics via operator learning for Alzheimer’s disease progression},
  author={Wang, Jindong and Mao, Yutong and Liu, Xiao and Hao, Wenrui and Alzheimer’s Disease Neuroimaging Initiative},
  journal={npj Systems Biology and Applications},
  year={2026},
  publisher={Nature Publishing Group UK London}
}

@article{zheng2022data,
  title={Data-driven causal model discovery and personalized prediction in Alzheimer's disease},
  author={Zheng, Haoyang and Petrella, Jeffrey R and Doraiswamy, P Murali and Lin, Guang and Hao, Wenrui and Alzheimer’s Disease Neuroimaging Initiative},
  journal={NPJ digital medicine},
  volume={5},
  number={1},
  pages={137},
  year={2022},
  publisher={Nature Publishing Group UK London}
}

@article{ferrell2012bistability,
  title={Bistability, bifurcations, and Waddington's epigenetic landscape},
  author={Ferrell, James E},
  journal={Current biology},
  volume={22},
  number={11},
  pages={R458--R466},
  year={2012},
  publisher={Elsevier}
}

@article{feinberg2023epigenetics,
  title={Epigenetics as a mediator of plasticity in cancer},
  author={Feinberg, Andrew P and Levchenko, Andre},
  journal={Science},
  volume={379},
  number={6632},
  pages={eaaw3835},
  year={2023},
  publisher={American Association for the Advancement of Science}
}

@inproceedings{zhang2025learning,
  title={Learning stochastic dynamics from snapshots through regularized unbalanced optimal transport},
  author={Zhang, Zhenyi and Li, Tiejun and Zhou, Peijie},
  booktitle={International Conference on Learning Representations},
  volume={2025},
  pages={19888--19919},
  year={2025}
}

@article{wang2025mathematical,
  title={Mathematical modeling and solution landscape reveal cancer progression dynamics in tumor ecological microenvironment},
  author={Wang, Shun and Wang, Tengfei and Wu, Shuonan and Zhang, Lei and Zou, Xiufen},
  journal={SIAM Journal on Applied Mathematics},
  volume={85},
  number={1},
  pages={50--77},
  year={2025},
  publisher={SIAM}
}

@article{lang2021landscape,
  title={Landscape and kinetic path quantify critical transitions in epithelial-mesenchymal transition},
  author={Lang, Jintong and Nie, Qing and Li, Chunhe},
  journal={Biophysical Journal},
  volume={120},
  number={20},
  pages={4484--4500},
  year={2021},
  publisher={Elsevier}
}

@article{lv2014constructing,
  title={Constructing the energy landscape for genetic switching system driven by intrinsic noise},
  author={Lv, Cheng and Li, Xiaoguang and Li, Fangting and Li, Tiejun},
  journal={PLoS one},
  volume={9},
  number={2},
  pages={e88167},
  year={2014},
  publisher={Public Library of Science San Francisco, USA}
}

@article{yin2021searching,
  title={Searching the solution landscape by generalized high-index saddle dynamics},
  author={Yin, Jianyuan and Yu, Bing and Zhang, Lei},
  journal={Science China Mathematics},
  volume={64},
  number={8},
  pages={1801--1816},
  year={2021},
  publisher={Springer}
}

@article{weinan2002string,
  title={String method for the study of rare events},
  author={Weinan, E and Ren, Weiqing and Vanden-Eijnden, Eric},
  journal={Physical Review B},
  volume={66},
  number={5},
  pages={052301},
  year={2002},
  publisher={APS}
}

@article{cang2023screening,
  title={Screening cell--cell communication in spatial transcriptomics via collective optimal transport},
  author={Cang, Zixuan and Zhao, Yanxiang and Almet, Axel A and Stabell, Adam and Ramos, Raul and Plikus, Maksim V and Atwood, Scott X and Nie, Qing},
  journal={Nature methods},
  volume={20},
  number={2},
  pages={218--228},
  year={2023},
  publisher={Nature Publishing Group US New York}
}

@article{cang2025synchronized,
  title={Synchronized optimal transport for joint modeling of dynamics across multiple spaces},
  author={Cang, Zixuan and Zhao, Yanxiang},
  journal={SIAM journal on applied mathematics},
  volume={85},
  number={1},
  pages={341--365},
  year={2025},
  publisher={SIAM}
}

@article{xu2014energetic,
  title={An energetic variational approach for ion transport},
  author={Xu, Shixin and Sheng, Ping and Liu, Chun},
  journal={arXiv preprint arXiv:1408.4114},
  year={2014}
}

@article{gangbo2019unnormalized,
  title={Unnormalized optimal transport},
  author={Gangbo, Wilfrid and Li, Wuchen and Osher, Stanley and Puthawala, Michael},
  journal={Journal of Computational Physics},
  volume={399},
  pages={108940},
  year={2019},
  publisher={Elsevier}
}
\bibliographystyle{ws-m3as}

\newpage
\appendix
\section{Appendices}
\subsection{Hamilton-Jaccobi Equation of Dynamic Optimal Transport}
The Monge-Kantorovich problem (MKP) is also called as optimal transport problem.The following theorem demonstrates the existence of $M$ as \cite{benamou2002monge}:
\begin{theorem}
There is a unique optimal mapping $\bar{M}$ defined on the support of $p_1(x)$ satisfying Monge form. The mapping $M$ is characterized as the unique mapping from this class which can be written as the gradient of a convex potential $\Phi$:
\begin{equation}
\bar{M}(x) = \nabla \Phi(x).
\end{equation}
\label{th1}
\end{theorem}
Using the Theorem \ref{th1} and linear interpolation to time $t$, given velocity field $u(x,t)$, the optimal flow map for MKP is established as:
\begin{equation}
x(X,t) = X + \frac{t}{T}(\nabla\Phi(X)-X), u(x,t) = \frac{1}{T}(\nabla\Phi(X)-X),
\end{equation}
and the Wasserstein distance of trajectory from $x(X,0)$ to $x(X,t)$ is expressed as
\begin{align}
W_2^2(\rho(x,0),\rho(x,T))&= \int_{\Omega} \|x(X,0) - x(X,T)\|^2  \rho(X,0)dX, \nonumber \\
&=T^2\int_{\Omega} \|u(x,t)\|^2  \rho(X,0)dX \nonumber \\
&=T\int_0^T\int_{\Omega} \|u(x,t)\|^2  \rho(x,t)dxdt \end{align}
By integrating the flow map and Wasserstein distance with Monge form, Benamou-Brenier (BB) formula \cite{benamou2000computational} is introduced that the Monge form of Wasserstein distance is equal to the optimization problem (Eq.\ref{eq:2}).

We introduce a scalar Lagrange multiplier 
$\phi(x,t)$, a scalar potential function, to impose the constraint as
\begin{equation}
\mathcal{L}[\rho,u,\phi] = T \int_0^T \int_{\Omega} \rho \|u\|^2 + \phi(x,t)(\partial_t \rho + \nabla \cdot (\rho u)) \, dx \, dt.
\end{equation}
Then, partial integration is performed separately over time and space as:
\begin{equation}
\int_0^T \phi\partial_t \rho dt = \phi(x,T)\rho(x,T)-\phi(x,0)\rho(x,0)-\int_0^T \rho\partial_t \phi dt,
\end{equation}
\begin{equation}
\int_{\Omega} \phi\nabla \cdot (\rho u)dx = -\int_{\Omega}\rho u\cdot \nabla\phi dx .
\end{equation}
Herein, the spatial boundary is assumed to be $\rho|_{\partial\Omega}=0$. The Lagrangian $\mathcal{L}[\rho,u,\phi]$ becomes:
\begin{equation}
\mathcal{L}[\rho,u,\phi] = T \int_0^T \int_{\Omega} \rho( \|u\|^2-\partial_t \phi- u\cdot \nabla\phi) dxdt + T \int_{\Omega} [\phi(x,T)\rho(x,T)-\phi(x,0)\rho(x,0)]dx,
\end{equation}
and the square in $u$ is completed as
\begin{align}
\mathcal{L}[\rho,u,\phi] =& T \int_0^T \int_{\Omega} \rho[(u-\frac{1}{2}\nabla\phi)^2-\partial_t \phi - \frac{1}{4}|\nabla\phi|^2) dxdt \nonumber\\
&+ T \int_{\Omega} [\phi(x,T)\rho(x,T)-\phi(x,0)\rho(x,0)]dx.
\end{align}
The minimization to Lagrangian to $u$ and $\rho$ is achieved at $u=\frac{1}{2}\nabla\phi$, and the variation of Lagrangian to $\rho$ as 
\begin{equation}
\frac{\delta\mathcal{L}}{\delta \rho} = \partial_t \phi + \frac{1}{4}|\nabla\phi|^2 =0.
\end{equation}
So, the dual Lagrangian is presented as
\begin{equation}
\sup_{\phi} T \int_{\Omega} [\phi(x,T)\rho(x,T)-\phi(x,0)\rho(x,0)]dx
\end{equation}
\[
\text{s.t.} \quad
 \partial_t \phi + \frac{1}{4}|\nabla\phi|^2 =0.
\]
If the time $t$ is rescaled as $2t$, the Hamilton-Jacobi equation to $\phi$ is obtained as 
\begin{equation}
\quad \partial_t \phi + \frac{1} {2}|\nabla\phi|^2 =0.
\end{equation}

\subsection{Tables}
\begin{table}[ht]
\centering
\label{tab:1}
\caption{Simulation parameters for the three gene regulatory motif.}
\label{tab:gene_params}
\begin{tabular}{lcl}
\toprule
\textbf{Parameter} & \textbf{Value} & \textbf{Description} \\
\midrule
$\alpha_1$    & $0.5$              & Self-activation strength for $x_1$. \\
$\gamma_1$    & $0.5$              & Inhibition strength exerted by $x_3$ on $x_1$. \\
$\alpha_2$    & $1$                & Self-activation strength for $x_2$. \\
$\gamma_2$    & $1$                & Inhibition strength exerted by $x_3$ on $x_2$. \\
$\alpha_3$    & $1$                & Self-activation strength for $x_3$. \\
$\gamma_3$    & $10$               & Half-saturation constant in the inhibition term. \\
$\delta_1$    & $0.4$              & Degradation rate for $x_1$. \\
$\delta_2$    & $0.4$              & Degradation rate for $x_2$. \\
$\delta_3$    & $0.4$              & Degradation rate for $x_3$. \\
$\eta_1$      & $0.05$             & Noise intensity for $x_1$. \\
$\eta_2$      & $0.05$             & Noise intensity for $x_2$. \\
$\eta_3$      & $0.01$             & Noise intensity for $x_3$. \\
$\eta_d$      & $0.014$            & Noise intensity for perturbations during cell division. \\
$\beta$       & $1$                & External signal activating $x_1$ and $x_2$. \\
\bottomrule
\end{tabular}
\end{table}
\newpage
\subsection{EVMDOT Algorithm to the Fokker--Planck equation}
\begin{algorithm}[ht]
\caption{Applying EVMDOT to the Fokker--Planck equation}
\label{alg:1}
\begin{algorithmic}[1]
\renewcommand{\algorithmicrequire}{\textbf{Input:}}
\renewcommand{\algorithmicensure}{\textbf{Output:}}
\Require Density trajectories $\{\rho^{(m)}(\mathbf{x},t_k)\}_{k=1}^{K}$ for $m=1,\dots,M$;
         noise intensity $\sigma$; domain $\Omega=\prod_{j=1}^{d}[a_j,b_j]$ with grid sizes $\{N_j\}$;
         time nodes $\{t_k\}_{k=1}^{K}$; trajectory weights $\{\lambda^{(m)}\}$; learning rate $\eta$; epochs $N$
\Ensure  Trained potential network $\psi_{nn}(x;\theta_\psi^{*})$
\State Initialize network parameters $\theta_\psi$
\State \(\triangleright\)~\textit{Quadrature Setup (Composite Trapezoidal Rule)}
\State Build $1$D weights $w_{j,i_j}$ by Eq.~\eqref{eq:9} for each direction $j$
\State Form tensor-product weights $W_{\mathbf{h}}=\prod_{j=1}^{d} w_{j,i_j}$ on grid points $\mathbf{x}_{\mathbf{h}}$
\State \(\triangleright\)~\textit{Energy-Dissipation Training}
\For{$\text{epoch}=1$ \textbf{to} $E$}
    \For{$m=1$ \textbf{to} $M$}
        \For{$k=1$ \textbf{to} $K$}
            \State Evaluate $\psi_{nn}({x}_{h};\theta_\psi)$ and $\nabla\psi_{nn}({x}_{h};\theta_\psi)$ on the grid
            \State Compute free energy $F^{(m)}(t_k)$ by Eq.~\eqref{eq:7}
            \State Compute dissipation $D^{(m)}(t_k)$ by Eq.~\eqref{eq:7}
        \EndFor
        \State Assemble discrete trajectory loss $\mathcal{L}^{(m)}$ by Eq.~\eqref{eq:8}
        \Statex \hspace{\algorithmicindent}\hspace{\algorithmicindent}%
        $\displaystyle \mathcal{L}^{(m)}=\Big(F^{(m)}(t_K)-F^{(m)}(t_1)
        +\sum_{k=1}^{K-1}\tfrac{\sigma^2}{4}\,\Delta t_k\big(D^{(m)}(t_k)+D^{(m)}(t_{k+1})\big)\Big)^{2}$
    \EndFor
    \State Compute total loss $\mathcal{L}=\sum_{m=1}^{M}\lambda^{(m)}\mathcal{L}^{(m)}$
           \Comment{Eq.~\eqref{eq:8}}
    \State Update $\theta_\psi \gets \theta_\psi-\eta\,\nabla_{\theta_\psi}\mathcal{L}$
\EndFor
\State \Return $\theta_\psi^{*}$, i.e. the learned potential $\psi_{nn}(x;\theta_\psi^{*})$
\end{algorithmic}
\end{algorithm}

\subsection{Hyperparameters on neural network}
For all examples, we use a constant weighting function $\lambda^{(m)} \equiv 1$ in the loss functions. The potential is parameterized by a fully connected neural network with three hidden layers of $32$ nodes each, using the GeLU activation function. The network is trained with the AdamW optimizer at a learning rate of $3 \times 10^{-4}$.

\subsection{Generation of Simulation Density Data}
We generate the simulation data by numerically solving the Fokker--Planck equation (Eq.~\eqref{eq:10}) on a uniform Cartesian grid using a finite-volume scheme. Time integration is carried out with an explicit forward Euler method, while the convection (drift) term is discretized using the Donor-Cell upwind flux and the diffusion term with standard second-order central differences \cite{ames2014numerical}. The explicit time step is chosen to satisfy the Courant--Friedrichs--Lewy (CFL) stability condition \cite{ames2014numerical}, ensuring a stable and non-negative evolution of the probability density.

\begin{algorithm}[htbp]
\caption{Finite-volume simulation of the Fokker--Planck equation}
\label{alg:fp_solver}
\begin{algorithmic}[1]
\Require initial density $\rho^{0}$, drift $a(x)$, noise intensity $\sigma$, grid spacing $h$, time step $\Delta t$, number of steps $N_t$, dimension $d\in\{2,3\}$
\Ensure time series $\{\rho^{n}\}_{n=0}^{N_t}$
\State initialize $\rho^{0}$ on the uniform Cartesian grid
\State verify the CFL condition $\Delta t \le \min\!\Big(\dfrac{h}{\max_x \|a(x)\|_1},\ \dfrac{h^{2}}{d\,\sigma^{2}}\Big)$
\For{$n = 0,1,\ldots,N_t-1$}
  \For{each interior grid point $\mathbf{i}$}
    \State \textbf{convection:} for each direction, compute the Donor-Cell upwind fluxes
    \Statex \hspace{1.5em} $F_{+}=\max(a_{+},0)\,\rho_{\mathbf{i}}+\min(a_{+},0)\,\rho_{\mathbf{i}+\mathbf{e}}$, \quad
    $F_{-}=\max(a_{-},0)\,\rho_{\mathbf{i}-\mathbf{e}}+\min(a_{-},0)\,\rho_{\mathbf{i}}$
    \State \textbf{diffusion:} compute the central-difference Laplacian $(\Delta_h\rho)_{\mathbf{i}}$ ($5$-point if $d=2$, $7$-point if $d=3$)
    \State \textbf{update:} $\rho_{\mathbf{i}}^{n+1} \gets \rho_{\mathbf{i}}^{n} - \dfrac{\Delta t}{h}\sum_{\text{dir}}(F_{+}-F_{-}) + \dfrac{\sigma^{2}\Delta t}{2h^{2}}(\Delta_h\rho)_{\mathbf{i}}^{n}$
  \EndFor
  \State \textbf{boundary:} set ghost nodes by even reflection,
$\rho_{\text{ghost}}^{n}=\rho_{\text{adjacent interior}}^{n}$, equivalently set
all boundary-face fluxes to zero (Eqs.~\eqref{eq:ghost_2d}--\eqref{eq:zeroflux}),
enforcing $\nabla\rho\cdot\mathbf{n}|_{\partial\Omega}=0$
\EndFor
\end{algorithmic}
\end{algorithm}

\paragraph{Two-dimensional case.}
On a uniform grid $x_{i}=x_0+i\,h$, $y_{j}=y_0+j\,h$ with spacing $h$, we use a five-point stencil in which each node $(i,j)$ is updated from its four nearest neighbors together with itself. Writing $\rho_{i,j}^{n}\approx\rho(x_i,y_j,t_n)$ and denoting the drift by $a=(a^x,a^y)$, the update reads
\begin{equation}
\rho_{i,j}^{n+1}=\rho_{i,j}^{n}
-\frac{\Delta t}{h}\Big(F^{x}_{i+\frac12,j}-F^{x}_{i-\frac12,j}+F^{y}_{i,j+\frac12}-F^{y}_{i,j-\frac12}\Big)
+\frac{\sigma^2\,\Delta t}{2\,h^2}\,\big(\Delta_h \rho\big)_{i,j}^{n},
\end{equation}
where $\big(\Delta_h\rho\big)_{i,j}^{n}=\rho_{i+1,j}^{n}+\rho_{i-1,j}^{n}+\rho_{i,j+1}^{n}+\rho_{i,j-1}^{n}-4\rho_{i,j}^{n}$ is the five-point Laplacian, and the Donor-Cell upwind flux across the interface $(i+\tfrac12,j)$ is
\begin{equation}
F^{x}_{i+\frac12,j}=\max\!\big(a^{x}_{i+\frac12,j},0\big)\,\rho_{i,j}^{n}
+\min\!\big(a^{x}_{i+\frac12,j},0\big)\,\rho_{i+1,j}^{n},
\end{equation}
with $F^{y}$ defined analogously in the $y$-direction.

\paragraph{Three-dimensional case.}
On a uniform grid with the additional coordinate $z_{k}=z_0+k\,h$, the five-point stencil is replaced by a seven-point stencil, in which each node $(i,j,k)$ is updated from its six nearest neighbors together with itself. The update takes the same form,
\begin{equation}
\rho_{i,j,k}^{n+1}=\rho_{i,j,k}^{n}
-\frac{\Delta t}{h}\sum_{d\in\{x,y,z\}}\!\Big(F^{d}_{+}-F^{d}_{-}\Big)
+\frac{\sigma^2\,\Delta t}{2\,h^2}\,\big(\Delta_h \rho\big)_{i,j,k}^{n},
\end{equation}
where $\big(\Delta_h\rho\big)_{i,j,k}^{n}$ is now the seven-point Laplacian summing the six nearest neighbors minus $6\rho_{i,j,k}^{n}$, and the Donor-Cell upwind fluxes $F^{d}_{\pm}$ in each direction $d$ are constructed as in the two-dimensional case.

\paragraph{No-flux boundary discretization.}
The no-flux boundary condition $\nabla\rho\cdot\mathbf{n}|_{\partial\Omega}=0$ is
imposed through the ghost-cell (mirror) technique, which guarantees that no
probability mass leaves the domain and preserves the total mass
$\int_\Omega\rho\,dx$ exactly at the discrete level. For each boundary face we
introduce a layer of ghost nodes outside $\Omega$ and assign them by even
reflection of the adjacent interior values, so that the normal derivative
vanishes to second order at the boundary. In the two-dimensional case, for the
$x$-boundaries this reads
\begin{equation}
\label{eq:ghost_2d}
\rho_{-1,j}^{n}=\rho_{0,j}^{n},
\qquad
\rho_{N_x+1,j}^{n}=\rho_{N_x,j}^{n},
\end{equation}
and analogously $\rho_{i,-1}^{n}=\rho_{i,0}^{n}$,
$\rho_{i,N_y+1}^{n}=\rho_{i,N_y}^{n}$ for the $y$-boundaries; the
three-dimensional case adds the same reflection in $z$. Equivalently, the
numerical flux through every boundary face is set to zero,
\begin{equation}
\label{eq:zeroflux}
F^{x}_{-\frac12,j}=F^{x}_{N_x+\frac12,j}=0,
\qquad
F^{y}_{i,-\frac12}=F^{y}_{i,N_y+\frac12}=0,
\end{equation}
so that both the upwind convection flux and the central-difference diffusion flux
vanish normal to $\partial\Omega$. With the reflected ghost values
\eqref{eq:ghost_2d}, the discrete Laplacian at a boundary node automatically
drops the outward neighbor; for instance, at the left $x$-boundary $(0,j)$,
\begin{equation}
\big(\Delta_h\rho\big)_{0,j}^{n}
=\rho_{1,j}^{n}+\rho_{0,j}^{n}+\rho_{0,j+1}^{n}+\rho_{0,j-1}^{n}-4\rho_{0,j}^{n}
=\rho_{1,j}^{n}+\rho_{0,j+1}^{n}+\rho_{0,j-1}^{n}-3\rho_{0,j}^{n},
\end{equation}
recovering the homogeneous Neumann (no-flux) stencil.
\newpage
\subsection{Supplementary Figures}
\begin{figure}[ht]
\centering
\includegraphics[width=1.0\linewidth]{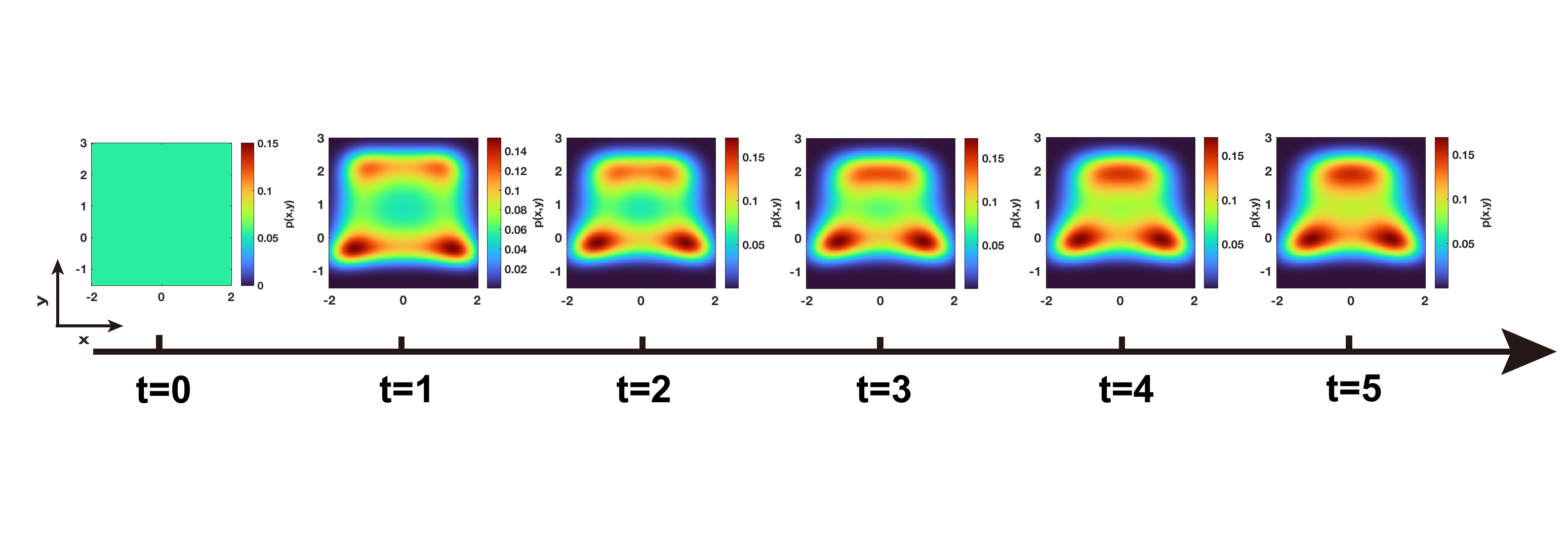}
\caption{\scriptsize\textbf{ Uniform Sampling in Early Phase.}
Simulation data generated from uniform initial distribution using the Fokker--Planck equation (Eq.~\eqref{eq:10}) together with the benchmark model (Eq.~\eqref{eq:11}); the colorbar indicates the probability density $\rho$. The value of $\sigma^2$ is set to 0.2.}
\label{fig:S3}
\end{figure}


\end{document}